\documentclass[12pt]{article}
\usepackage{amsfonts}
\usepackage{amssymb}
\usepackage{xltabular}
\usepackage{mathrsfs}
\usepackage{srcltx}
\usepackage{amsmath}
\usepackage{amsthm}
\usepackage{amstext}
\usepackage{amsopn}
\usepackage{graphicx}
\usepackage{url}
\usepackage{color}
\usepackage{bm}
\usepackage{physics}
\usepackage{appendix}
\usepackage{lipsum}
\usepackage{float}
 \usepackage{hyperref}
\newtheorem{theorem}{Theorem}[section]
\newtheorem{lemma}[theorem]{Lemma}
\newtheorem{proposition}[theorem]{Proposition}
\newtheorem{corollary}[theorem]{Corollary}
\theoremstyle{definition}

\usepackage{placeins}
\theoremstyle{remark}

\usepackage{array, makecell} 
\numberwithin{equation}{section}
\newcommand{\ba}{\begin{array}}
\newcommand{\ea}{\end{array}}
\newcommand{\f}{\frac}

\newcommand{\ds}{\displaystyle}
\newcommand{\g}{\gamma}
\usepackage[affil-it]{authblk} 
\usepackage{etoolbox}
\usepackage{lmodern}

\makeatletter
\patchcmd{\@maketitle}{\LARGE \@title}{\fontsize{16}{19.2}\selectfont\@title}{}{}
\makeatother

\usepackage{comment}

  {
      \theoremstyle{plain}
      \newtheorem{assumption}{Assumption}
  }

\begin{document}
\date{}
\title{ An Epidemic Compartment Model for Economic Policy Directions for Managing Future Pandemic }

\author[1]{Zachariah Sinkala\thanks{corresponding author : Zachariah.Sinkala@mtsu.edu } }
\author[1]{Vajira Manathunga\thanks{vmanathunga@mtsu.edu}}
\author[2]{Bichaka Fayissa \thanks{bichaka.fayissa@mtsu.edu} }
\affil[1]{Department of Mathematical Sciences, Middle Tennessee State University}
\affil[2]{Department of Economics and Finance, Middle Tennessee State University}
\maketitle

\begin{abstract}
\begin{footnotesize}
 In this research, we develop a framework to analyze the interaction between the economy and the Covid-19 pandemic using an extension of SIR epidemic model. At the outset, we assume there are two health related investments including general medical expenditures and the other for a direct investment for controlling the pandemic. We incorporate the learning dynamics associated with the management of the virus into our model. Given that the labor force in a society depends on the state of the epidemic, we allow birth, death, and vaccination to occur in our model and assume labor force consists of the susceptible, vaccinated, and recovered individuals.  We also assume parameters in our epidemic compartmental model depend on investment amount for directly controlling the epidemic, the health stock of individual representative agents in the society, and the knowledge or learning about the epidemic in the community. By controlling consumption, the general medical expenditure, and the direct investment of funds for controlling the epidemic, we optimize the utility realized by the representative individuals because of consumption.  This problem is nontrivial since the disease dynamics results in a non-convex optimization problem. 
\end{footnotesize}
\end{abstract}
\begin{footnotesize}
\textbf{\indent Keywords}: SIR Model, Epidemiology,  Existence of equilibrium, Steady States, Economic models

\indent \textbf{AMS Classification (MSC2020): 91B15, 91B16, 91B50, 91B55 } 
\end{footnotesize}

\section{Introduction}
 COVID 19 is the most significant event for humans in this century. Compare to the impact of  two world wars and several epidemics occurred in last century, this has equivalent or may be more impact on the society. The impact can be social, economical or both. Covid 19 impact on economies around the world is significant. For example, \cite{maital2020global} summarized recent research and report on the global economic impact of the COVID 19. However, the severity of the impact changed from one region to the other. Because of that, a slew of research articles focused on economic impact based on the region or country was published. For example, \cite{kumar2020social} discuss the social impact of COVID 19 in India, \cite{abiad2020economic} discuss the impact on developing Asia,  \cite{lone2020covid} put emphasis on African region , \cite{suryahadi2020impact} discuss the COVID 19 impact on poverty in Indonesia, and \cite{hu2020intersecting} research on COVID 19 impact on UK economy name to few. 
 
 Many countries opt to allocate money to directly combat the virus during this time and use many mitigation policies. These prevention, intervention, and mitigation policies and their impacts were investigated in several articles. A review of current interventions were done in \cite{pradhan2020review}, effects of interventions were discussed in \cite{soltesz2020effect}, A model for COVID 19 intervention policies in Italy were given in \cite{giordano2020modelling}, A model to simulate the health and economic effect of social distancing given in \cite{silva2020covid}. One of tragedy in COVID was the loss of human lives. Labor-wise, the workforce shrank significantly due to death, infection, and movement restrictions, which significantly affect the growth of the economy. We have mentioned that COVID has direct impact on labor or human capital in above. This indirectly cause the shortage in physical capital. For example current ``chips shortage" \cite{wu2021analysis}, and port blockage \cite {notteboom2021iaph} discuss implication of COVID on physical capital. According to the World Bank \cite{maliszewska2020potential}, baseline GDP of the world would fall below   2  percent from pre-pandemic baseline, for developing countries by  2.5 percent,  and  for industrial countries  by  1.8  percent.
 
History has shown that humans eventually learn how to either eradicate or live with diseases. For example, smallpox, rinderpest are completely eradicated and we learned to live with AIDS using prevention methods. Similar to this, learning occurred over the time about COVID 19. At the early pandemic stage, countries opted to close schools, workplaces, and economies; however, as time passed, these economies decided to open again and learned to live with the COVID. Even though new variants such as Delta and Omicron are becoming more prevalent, many countries opt to keep open, using the knowledge accumulated over the time to handle these new variants. Therefore it is clear that we should look into learning about the disease that occurs during this time. Learning does not have to be strictly pharmaceutical, but it can be nonpharmaceutical learning, such as wearing masks, keeping 6 feet apart from others, etc.  Also, another aspect of the COVID is the introduction of the vaccine. Even though the vaccine was commonly available in developed countries, subgroups started to show resistance to getting vaccinated. Thus, vaccine hesitancy is becoming a key issue in these countries, which implies COVID 19 may become endemic over time.

Infectious diseases may bring shocking effects to economies. Therefore understanding the impact of epidemics on economies is crucial. One approach would be to model the economy's dynamics without disease and introduce disease as an exogenous shock to the model. However, this would prevent us from understanding the two-way interactions between the disease dynamics and the economy. For example, as we mentioned earlier, when the government directly allocates money to combat a deadly disease, it impacts disease transmission and other aspects of the disease. On the other hand incidence of disease negatively affect the labor force, which causes lower physical output in the physical production function.  Therefore, we believe it is best to model the disease and the economic dynamics together.  In this research, we model the accumulation of physical and health capital through the neo-classical growth model.  COVID 19 brought renewed attention to compartmental models and their applications.We endogenize disease and use SIRV (Susceptible, Infected, Recovered, Vaccinated) model for disease dynamics.The interaction between the neo-classical growth model and the SIRV model is introduced as follows: First, we assume the labor force consists of people who are not infected. Second, we consider two parameters in the SIRV model, namely transmission rate and recovery rate; both depend on the money allocated  to combat the disease, the knowledge in the society about how to control the disease, and health capital. Thus, disease and economic dynamics are needed to solve simultaneously. Another aspect of our model is dynamics of knowledge for controlling the disease.

We have mentioned that we incorporate learning dynamics into our model.    The learning or knowledge production function uses direct investment to control the epidemic and existing knowledge to create more understanding about managing the outbreak. Therefore the learning by controlling function help to reduce the cost associated with disease control and improve the health stock. The process is intuitive in many ways. Any society which faces an outbreak eventually accumulates or created knowledge about how to control it, which will be used in the next epidemic. 

In our model, we consider two types of health expenditure. The health capital production function produces health. In order to produce health stock, the health production function uses general health expenditure, which results in health, or more healthy people to be accurate. Healthy people are less prone to become infected. Therefore this health is required. The second form of health expenditure comes from direct allocation to combat the deadly disease, such as money to do research and find vaccines and other cures. This is also crucial when battling a virus-like COVID 19.  

Several previous studies focused on modeling disease and economic dynamics together.  The effect of recurring disease on an economy that follows SIS dynamics was investigated in \cite{goenka2014infectious}. In that paper, the authors tried to understand the best society can do to control the disease transmission while optimizing the utility gain by the consumption. The authors showed that a steady-state with zero health expenditure could be an optimal solution to the problem. However, the authors never considered a direct investment to control the epidemic or learning that occurs when controlling the disease. Also SIS model assumes no immunity from the disease compared to SIRV model we used. 

SIR model with two way interaction between the economy and the disease was investigated in \cite{boris}. In this model, the authors introduced learning by controlling and direct investment to control the epidemic. However, the economy consists of a single sector with input as health capital. We have a two-sector economy where physical capital and health capital are both produced compared to this model. Also, the authors assumed that death could occur only when in infected state, but we assumed death could occur at any compartment. The paper lacks rigorous mathematical proof, which establishes necessary and sufficient conditions of optimal solutions and calculation of steady states.

The authors in \cite{d2009optimal} used the SI epidemic model and claimed that a minimal level of labor is necessary to reduce the long-term prevalence rate of the epidemic. If this minimal labor is not reached, then prevention would be temporary and may not be optimal to take. On the other hand, if enough labor exists, allocating resources to prevent the disease is feasible but not optimal. 

Epidemic models and their interactions with the economy was studied in \cite{goenka2020infectious, la2020optimal,d2021optimal}
More recently, SIR model, disease-related mortality and effect on economy was studied in  \cite{goenka2021sir}.

\section{Model}
\subsection{Epidemiology model}
\begin{figure}[h!]
\includegraphics[scale=0.6]{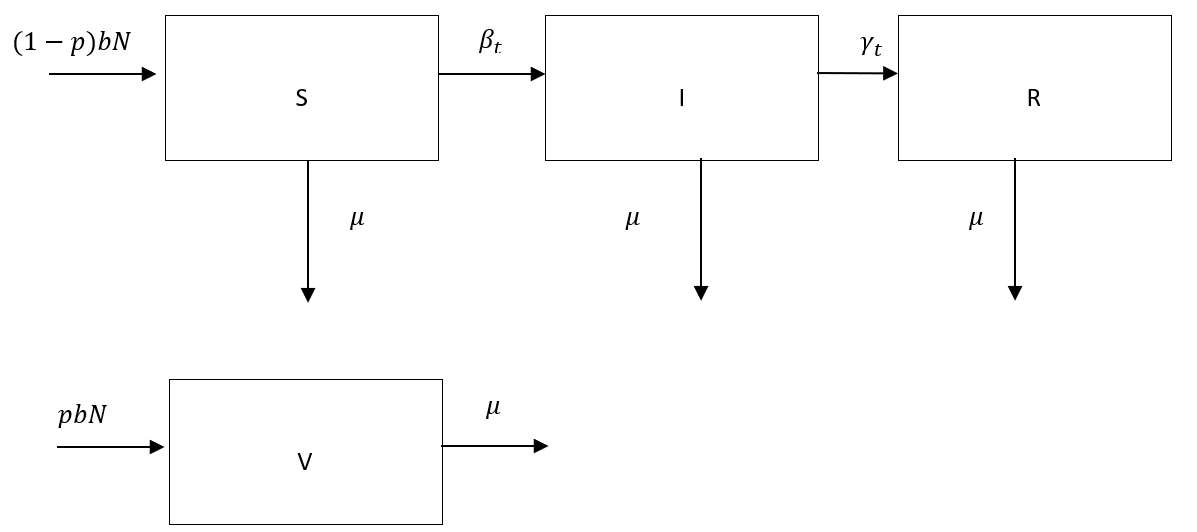}
\caption{{ \textbf{The model.} Schematic diagram representing the interactions among different stages of infection in the mathematical model SIRV: $S$ susceptible(uninfected), $I$, Infected, $R$, removed, $V$, Vaccinated}}
\label{sirdfig}
\end{figure}

We use modified SIR model with the vaccination to describe the dynamics of the disease. The model is given in figure \ref{sirdfig}.

We assume under epidemic outbreak, direct investment, $A_t$, to control epidemic will occur and some existing knowledge of controlling of the disease, $e_0$. We assume,  as time pass by learning about disease control $e_t$, will occur where $e_t$ is the accumulated knowledge of disease\cite{boris}.    


We denote the knowledge production function: learning by controlling as $E(A_t, e_t)$. We also assume knowledge about disease control can be depreciated over time \cite{knowledgedep}. Thus law of motion for knowledge production is given by
\begin{equation}\label{learn}
\dv{e_t}{t}=E(A_t,e_t)-\delta_E e_t .
\end{equation}
In equation (\ref{learn}), function $E(A_t,e_t)$ is increasing in control measures $A_t$, existing knowledge about disease control $e_t$ and existing knowledge depreciated at constant rate $\delta_E\in (0,1)$. 
We also assume disease transmission rate $\beta(A_t,e_t,h_t)$, recovery rate and $\gamma(A_t,e_t,h_t)$ are depend on $A_t, e_t$ and $h_t$ where $h_t$ is the health capital. In simplified assumption, we allow the transmission rate $\beta$ to depend on experience in controlling the disease. In our model (\ref{sirdraw}), we use mass action incidence principle. { We denote natural birth rate by $b$, natural death rate as $\mu$ and  assume these are constants}. { We also assume death due to the disease is negligible compare to the population size}. As an example this model can be used in  countries where epidemic death is less compared to total population. Assume that vaccination is carried out at recruitment. Thus, fraction $0\le p\le 1$ of population is vaccinated at time $t$ and become immune to the disease, and only $(1-p)$ fraction of population enter into the susceptible class. Let $N_t$ denote the total population at time $t$.  Then, we have $S_t+I_t+R_t+V_t=N_t$. The respective dynamics of modified Susceptible-Infected-Recovered (SIR) epidemiology model with vaccination is  given by:
\begin{align} \label{sirdraw}
\begin{split}
    \dot{S}_t=\dv{S_t}{t} &=(1-p)bN_t-\mu S_t-\frac{\beta(A_t,e_t,h_t) S_tI_t}{N_t}\\ 
    \dot{I}_t=\dv{I_t}{t} &=\frac{\beta(A_t,e_t,h_t) S_t I_t}{N_t}-(\gamma(A_t,e_t,h_t)+\mu) I_t\\
    \dot{R}_t=\dv{R_t}{t} &= \gamma(A_t,e_t,h_t) I_t -\mu R_t\\
      \dot{V}_t=\dv{V_t}{t} &= pbN_t  -\mu V_t\\
    \dot{N}_t=\dv{N_t}{t} &= (b-\mu)N_t
    \end{split}
\end{align}
Define $s_t=\frac{S_t}{N_t}, i_t=\frac{I_t}{N_t}, r_t=\frac{R_t}{N_t} $ and $v_t=\frac{V_t}{N_t}$. Observe that 
\begin{equation}
\begin{split}
\dv{s_t}{t} &=\dv{\big(\frac{S_t}{N_t}\big)}{t} \\
&=\frac{\dot{S_t}N_t-S_t\dot{N_t}}{N_t^2}\\
&=\frac{\dot{S_t}}{N_t}-s_t\frac{\dot{N}_t}{N_t}\\
&=\bigg((1-p)b-\mu s_t -\beta(A_t,e_t,h_t)s_ti_t\bigg) -\bigg((b-\mu)s_t \bigg)\\
&=(1-p)b-bs_t-(\beta(A_t,e_t,h_t))s_ti_t
\end{split}
\end{equation}
Similarly, we get dynamics of fractions of SIR epidemiology model  as follows:
\begin{align} \label{sird}
\begin{split}
    \dv{s_t}{t} &=(1-p)b-bs_t-(\beta(A_t,e_t,h_t))s_ti_t\\ 
    \dv{i_t}{t} &=\beta(A_t,e_t,h_t) s_t i_t-(\gamma(A_t,e_t,h_t)+b) i_t\\
    \dv{r_t}{t} &= \gamma(A_t,e_t,h_t) i_t -b r_t\\
    \dv{v_t}{t} &=pb-bv_t
    \end{split}
\end{align}

%
From this point onward, we would suppress $\beta(A_t,e_t,h_t), \gamma(A_t,e_t,h_t)$ and $\sigma(A_t,e_t,h_t)$ to be just $\beta,\gamma$ and $\sigma$ respectively. 
Before we continue to the economic model, we would like to analyse steady states of this epidemic model. The disease free equilibrium points \cite{martcheva}  are given by
\begin{equation} \label{dfree}
s^* = 1-p, i^* = 0, v^*=p, r^* = 0.
\end{equation}
 Jacobian of the disease free equilibrium for full system is given by
\[ 
\begin{bmatrix}
-b & -\beta(1-p) & 0 & 0 \\
0 & \beta(1-p) - b-\gamma & 0 & 0 \\
0 & \gamma & -b & 0 \\
0 & 0& 0 & -b
\end{bmatrix}
. \]

\noindent Eigenvalues of the disease free Jacobian are
\[\beta(1-p)-b-\gamma, -b, -b, -b .\]
The local stability of the solution depends on whether the eigenvalues of the Jacobian are negative or have negative real parts \cite{allen2007}. 
Observe that if $\beta(1-p) < (b+\gamma)$ then disease free steady state is stable otherwise unstable.
If fraction of the population vaccinated at recruitment, $p$, is zero, then reproduction number is $\mathcal{R}_0=\frac{\beta}{b+\gamma}$ \cite{martcheva}. Thus by vaccinating, we would reduce the original reproduction number by $(1-p)$. As mentioned in \cite{martcheva}, if we need the reproduction number to be less than 1, then at least $1-\frac{1}{\mathcal{R}_0}$ fraction of population need to be vaccinated otherwise the disease will invades the population. Next, endemic state equilibrium points \cite{martcheva} are given by
\begin{align}\label{endemic}
\begin{split}
s^* &= \frac{\gamma^*+b}{\beta^*}=\frac{1}{\mathcal{R}_0}\\
i^* &= \frac{(1-p)b}{\gamma^*+b}-\frac{b}{\beta^*} = \frac{(1-p)b}{\gamma^*+b}-\frac{b}{\mathcal{R}_0 (\gamma^*+b)}\\
v^* &= p \\
r^* &= 1-(s^*+i^*+v^*) . \\
\end{split}
\end{align}
We need $\mathcal{R}_0> 1$, and $(1-p)\mathcal{R}_0>  1$ to make sure $0\le s^*, i^*\le 1$ in endemic steady state. However former follow from the latter inequality. Jacobian of the endemic state equilibrium
\[
\begin{bmatrix}
\frac{b \beta p - b \beta}{b + \gamma} & -(\gamma+b) & 0 & 0 \\
-\frac{b \beta p + b^{2} - b \beta + b \gamma}{b + \gamma} & 0 & 0 & 0 \\
0 & \gamma & -b & 0 \\
0 & 0 & 0 & -b
\end{bmatrix}.
\]

\noindent Define, 
\footnotesize
$$x=\sqrt{b^{2} \beta^{2} p^{2} + 4 \, b^{4} - 4 \, b^{3} \beta + b^{2} \beta^{2} + 4 \, b \gamma^{3} + 4 \, {\left(3 \, b^{2} - b \beta\right)} \gamma^{2} + 4 \, {\left(3 \, b^{3} - 2 \, b^{2} \beta\right)} \gamma + 2 \, {\left(2 \, b^{3} \beta - b^{2} \beta^{2} + 4 \, b^{2} \beta \gamma + 2 \, b \beta \gamma^{2}\right)} p}$$
\normalsize
Then eigenvalues of this Jacobian given by $\frac{-(1-p)\beta b - x}{2 \, {\left(b + \gamma\right)}}, \frac{-(1-p)\beta b + x}{2 \, {\left(b + \gamma\right)}}, -b, -b$

\begin{proposition}
If $(1-p)\mathcal{R}_0>1$ then endemic state is stable.
\end{proposition}
\begin{proof}
Suppose $(1-p)\mathcal{R}_0=\frac{(1-p)\beta}{b+\gamma}>1$. Then, notice that all eigen values except $\frac{-(1-p)\beta b + x}{2 \, {\left(b + \gamma\right)}}$ are negative. Suppose this eigen value is nonnegative. Then the determinant of the Jacobian matrix should be nonpositive. Observe that determinant of the Jacobian of the endemic state is $$b^3((1-p)\beta-(b+\gamma))>0$$ Hence  when $\frac{(1-p)\beta}{b+\gamma}>1$, all eigen values are negative and endemic state is stable.
\end{proof}

Assume vaccination rate $p$ is constant and define $\mathcal{R}_{vac}=(1-p)\mathcal{R}_0$. Thus the system has one stable disease free steady state if $\mathcal{R}_{vac}<1$. Otherwise, the system has two steady states with stable endemic state and unstable disease free steady state. Hence, $\mathcal{R}_{vac}$ is the bifurcation point. Following \cite{martcheva}, we give following diagram to illustrate forward bifurcation.

\begin{proposition}\label{endemicsteady}
For given path of $A_t,e_t,h_t$ if $(1-p)\beta(A,e,h)<b+\gamma(A,e,h)$ then the system has one stable disease free steady state otherwise the system has two steady states with stable endemic state and unstable disease free steady state.
\end{proposition}

\begin{figure}[h]
\includegraphics[scale=0.65]{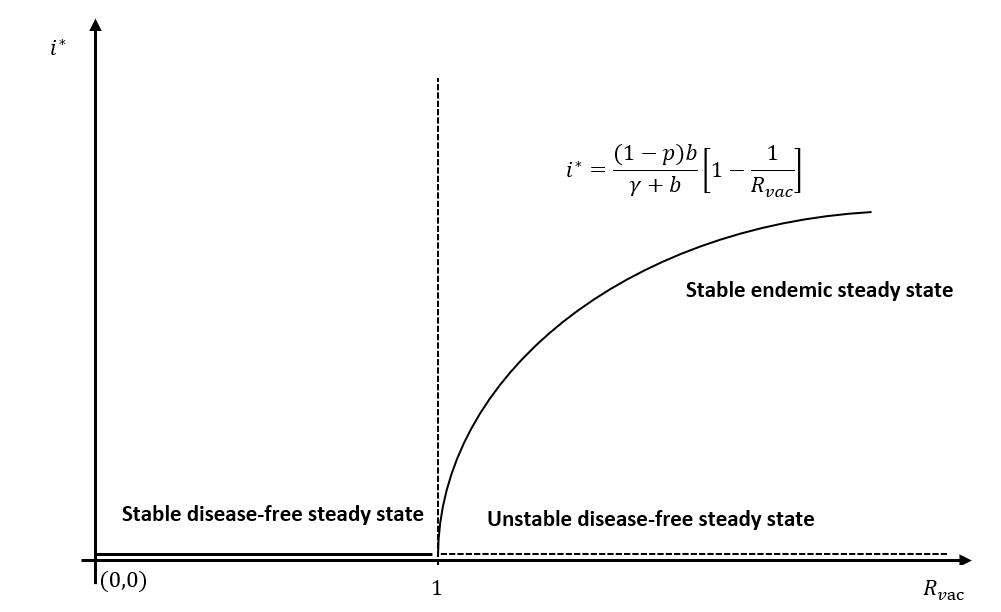}
\caption{Forward bifurcation diagram with respect to the reproduction number.
}
\label{bifurcation}
\end{figure}






\subsection{Economic model}
To avoid keeping track of individual health histories, we use the framework of a large representative consumer. We assume the economy is populated by non atomic identical consumers whose consumption is same irrespective to the health status and no need to keep track of the health record of individual consumers \cite{goenka2014infectious}. Thus, social planner can maximize the respective optimization problem for the representative consumer. We assume social planner try to maximize following utility function
\[\int_0^{\infty} u(c_t)e^{-\theta t}dt\]
where $\theta$ is the discount factor, $c_t$ is the consumption. Social planner face following problem: how to choose consumption $c_t$, general medical expenditure $m_t$ and direct epidemic control investment $A_t$ such that utility function is maximized. We assume general medical expenditure $m_t$ and direct epidemic control investment $A_t$ are mutually exclusive.
We assume there are three production functions: one for physical goods, one for health and one for creating knowledge controlling epidemic. Physical goods are either consumed, invested again, spend in health or depreciated. The production function $f$ of physical goods is depend on capital $k_t$, effective labor $l_t$. However under our model, effective labor 
\[l_t=s_t+v_t+r_t=1-i_t\]
Thus, we have 
\[\dv{k_t}{t}= f(k_t, l_t)-m_t-A_t-c_t-(\delta_K+b-\mu)k_t\]
where $\delta_K\in (0,1)$, denotes the depreciation of physical capital. The health capital production function $g$ depend on medical expenditure $m_t$ and hence we have
\[\dv{h_t}{t}=g(m_t)-(\delta_H+b-\mu)h_t\]
Here $\delta_H \in (0,1)$ denotes the depreciation of health capital.

\subsection{Assumptions}\label{assumptions}
We make following assumption regarding this model. 
\begin{assumption}Natural birth rate, disease death, and  natural death rates satisfy following inequality: $b\ge \mu$  \cite{d2008existence}.  \label{ass1}
    \end{assumption}
    \begin{assumption}
  Transmission rate function $\beta(A_t, h_t, e_t)$ satisfy following conditions: \label{ass2}
    \begin{enumerate}
        \item $\beta(\cdot):\mathbb{R}_{+}^3\to (0,1]$
        \item $\beta(0,0,0)=1$
        \item partial derivatives $\beta_1, \beta_2,\beta_3\le 0$ where $\beta_i$ denote partial derivative of $\beta(\cdot)$ respect to $i^{th}$ variable: This is justified since effective transmission rate should be decreased as $A_t, h_t, e_t$ increases.
        \item The second order partial derivatives, $\beta_{11},\beta_{22}, \beta_{33}\ge 0$. 
        \item $\beta\to \bar{\beta}$ as $(A_t,h_t,e_t)\to (0,0,0)$
    \end{enumerate}
    \end{assumption}
    \begin{assumption} Assuming infection last more than one day, recovery rate $\gamma(A_t, h_t, e_t)$ (the inverse of average infectious period) satisfy following conditions: \label{ass3}
       \begin{enumerate}
        \item $\gamma(\cdot):\mathbb{R}_{+}^3\to [0,1]$
        \item The first order partial derivatives $\gamma_1,\gamma_2,\gamma_3\ge 0$.
        \item The second order partial derivatives $\gamma_{11},\gamma_{22},\gamma_{33}\le 0$. 
        \item $\gamma\to \underline{\gamma}$ as $(A_t,h_t,e_t)\to (0,0,0)$
    \end{enumerate}
    \end{assumption}
    \begin{assumption} The physical production function $f(k_t, l_t):\mathbb{R}_+^2\to \mathbb{R}_+$ satisfy following Inada(Weak Inada) conditions \cite{goenka2014infectious}: \label{ass4}
    \begin{enumerate}
        \item $f\in C^{\infty}$ 
        \item The gradient of $f$ has positive components and the Hessian of $f$ is negative definite.
        \item $\lim_{k_t\rightarrow 0}f_1(k_t,l_t)=\infty,\,\,\lim_{k_t\rightarrow \infty}f_1(k_t,l_t)=0,\textup{ and } f(k_t,0)=f(0,l_t)=0, $ where $f_1$ denotes a partial derivative of with respect to $k_t$.
    \end{enumerate}
    \end{assumption}
    \begin{assumption}The health production function $g(m_t)$ satisfy the following Inada conditions\cite{inada1963two} :\\
    $g(\cdot):\mathbb{R}_{+}\rightarrow \mathbb{R}_{+} \textup{ is } C^{\infty},\textup{ with } \\g'>0,\,g''<0,\,\lim_{m_t\rightarrow 0}g'(m_t)<\infty, \lim_{m_t\to\infty}g'(m_t)=0 \textup{ and } g(0)=0$.
    \label{ass5}
    \end{assumption}
    \begin{assumption} The utility function $u(c_t):\mathbb{R}_+\to \mathbb{R}_+$  is $C^{\infty}$ with \\$u'>0,u''<0, \lim_{c_t\to 0}u'=\infty$. \cite{goenka2014infectious} \label{ass6}
    \end{assumption}
    \begin{assumption} There exist {  $\kappa_1$ such that $  \kappa_1 \in (0, \infty)$ and} $-\kappa_1 \le \ds\f{\dv{k_t}{t}}{k_t}$. \cite{goenka2014infectious} \label{ass7}
    \end{assumption}
    \begin{assumption}The knowledge production for disease control, $E(A_t,e_t):\mathbb{R}_+^2\to\mathbb{R}_+$ satisfy following conditions:
    \begin{enumerate}
        \item We assume direct investment to control the epidemic is finite (since resources are scarce). That mean, there exist $A>0$ such that $A_t\le A$ for all $t$. 
        \item $E(.,.) \in C^{\infty}$
         \item Since direct investment to control the epidemic, $A_t$ and the existing knowledge both are essential $E(A_t,0)=E(0,e_t)=0$
        \item $E_1>0,E_2>0, E_{11}<0, E_{22}<0$, $E_{12}=E_{21}$ and $E_{11}E_{22}-E_{12}E_{21}>0$, here $E_1, E_2$ denote partial derivative of $E(.,.)$ respect to $A_t$ and $e_t$.
        \item  $ \lim_{A_t\to \infty}E_1(A_t,e_t)=0$, $\lim_{e_t\to \infty}E_2(A_t,e_t)=0$  
    \item $E_2\to \bar{E}_2$ as $(A_t,e_t)\to (0,0)$
    \item $E_1\to \bar{E}_1$ as $(A_t,e_t)\to (0,0)$
    \end{enumerate} \label{ass8}
\end{assumption}
\subsection{Social Planner Problem}
Consider the following social planner problem($\mathcal{P}$):\\
Let $x_t=(c_t, m_t, A_t),y_t=(k_t, h_t, s_t, i_t,e_t)$ and assume the assumptions 1-9.
\begin{equation}\label{social}
    \max_{x_t,y_t}\int_0^{\infty} u(c_t)e^{-\theta t}dt 
\end{equation} 
subject to constraints\\
\begin{align}
    \dv{k_t}{t} &=f(k_t, l_t)-m_t-A_t-c_t-\bigg(\delta_K+b-\mu\bigg)k_t \label{eqk}\\
    \dv{h_t}{t} &=g(m_t)-\bigg(\delta_H+b-\mu\bigg)h_t \label{eqh}\\
    \dv{s_t}{t} &=(1-p)b-bs_t-\beta s_ti_t \label{eqs}\\
    \dv{i_t}{t} &= \beta s_t i_t-(\gamma+b)i_t\label{eqi}\\
    \dv{e_t}{t} &=E(A_t,e_t)-\delta_E e_t \label{eqe}\\
    l_t &=1-i_t\\
    k_t &>0,m_t\ge 0, h_t\geq 0,\,e_t>0,A_t\ge 0, s_t \ge 0, 1\ge l_t\\
    k_0 &>0,m_0\ge 0,h_0\geq 0,e_0\ge  0, A_0\geq 0,s_0>0, l_0>0 \label{eqineq}
\end{align}
where $x_t$ is the control vector variable, $y_t$ is the state vector variable.

\section{Existence of optimal solution}
In this section, we show that there exist optimal solution to the above mentioned social planer problem. Usual approach to show that existence of optimal solution is to use current Hamiltonian and then  apply either Mangasarian sufficient conditions \cite{mangasarian1966} or Arrow sufficient\cite{seierstad1986optimal} conditions. Under Mangasarian sufficient conditions, the current Hamiltonian must be concave w.r.t to all state  and control variables. This may not be true in our situation. For Arrow sufficient conditions, concavity with respect to the state variables of the Hamiltonian maximized with respect to control variable is required \cite{seierstad1977}. Again, this may not be true in our situation. For our problem, we rely on theorem 1 given in \cite{d2008existence}. First we restate the theorem 1 given in \cite{d2008existence} without proof and then show that how it can be applied to our problem.

Consider the following social planner problem($\mathcal{P}$)\cite{d2008existence}:\\
\[ \max_{x_t,y_t}\int_0^{\infty} u(x_t)e^{-\theta t} \,dt \]
subject to constraints\\
\begin{equation}\label{cong}
    \begin{cases}
    \textup{for all } t\geq 0,\, F(y_t) \leq \dv{y_t}{t} \leq G\left(\bar{y}_t,y_t,x_t\right)\\
    x_t\in \mathbb{R}_{+}^C,\, \bar{y}_t,y_t\in \mathbb{R}_{+}^{K}\\
    y_0\in \mathbb{R}_{+}^{K}, \theta >0 \textup{ are given},
    \end{cases}
\end{equation}
where $\bar{y_t}=(\bar{y}_1, \cdots,\bar{y}_K)$ is the external variable, $y_t=(y_1,\cdots,y_K)$ is the state variable,  $x_t=(x_1,\cdots,x_C)$ is the control variable and $\dv{y_t}{t}$ is the derivative of $y_t$. $C, K$  denote number of control variables and number of state variables in the optimization problem, which are finite.
Assume the following:
\begin{itemize}
\item[\textbf{A1}] The { function  $F$ is }  continuous on  $\mathbb{R}_{+}^K$ and function $G$ { is } continuous on  $\mathbb{R}_{+}^K\times\mathbb{R}_{+}^K\times \mathbb{R}_{+}^C $
\item[\textbf{A2}] For any $ j\in K,$ the function $ G^j$  is concave with respect to $x_t$.
\item[\textbf{A3}] There exists $( b_i\geq 0, B_i>0,\,i=1,2,...,C)$ and $( a_j\geq 0, A_j>0,\,j=1,2,...,K)$ such that:
\begin{enumerate}
    \item For all $t$, for all $j \quad \bar{y}_t^j\leq A_j e^{a_j t}$
    \item If $x_t$ and $y_t$ satisfy the following differential constraint for all $t$,
\[ F(y_t) \leq \dv{y_t}{t} \leq G\left(\bar{y}_t,y_t,x_t\right)  \] then for all $t$,
                   \begin{equation}
                    \begin{cases}
                     \text{ for all }  j,\quad y^j_t\leq A_j e^{a_j t}, \text{ and } |\dv{y_t^j}{t}|\leq A_j e^{a_j t} \\
                    \textup{ for all }   i,\quad x^i_t\leq B_i e^{b_i t}  \\
                    \end{cases}
                    \end{equation}

\end{enumerate}

\item[\textbf{A4}] $\theta >\sup\{b_i,i\in \{1,2,..,C\}\}$
\item[\textbf{A5}] The function $u$  is concave, non decreasing, and upper semi-continuous from 
        $\mathbb{R}_{+}^C$  into $\mathbb{R}\cup \{-\infty\}$ and finite valued on  $\mathbb{R}_{++}^C$
\end{itemize}

\begin{theorem}\label{exsistancemain}
Under assumptions $A1,\cdots, A5$, there exists a solution to the social planner's problem $(\mathcal{P})$ given in Equation (\ref{cong}) \cite{d2008existence}
\end{theorem}
We can apply above theorem \ref{exsistancemain}  to social planner problem given in Equation (\ref{social}). However, first we need to show that social planner problem (\ref{social}) satisfy assumption $A1,\cdots,A5$ 
\begin{lemma}\label{lemmaF}
Let $y_t=(k_t, h_t, s_t, i_t,e_t)$ be the state variables given in social planner problem in (\ref{social}). Then, there exists continuous function $F(\cdot):\mathbb{R}^5_+\to \mathbb{R}^5$ such that 
\[ F(y_t)\le \dot{y}_t\]
\end{lemma} 
\begin{proof} First we claim that there exists constants $\kappa_j>0$ such that $-\kappa_j y^j_t\le \dot{y}^j_t$. From section \ref{assumptions},  Assumption (\ref{ass7}), we know there exists constant $\kappa_1>0$ such that $-\kappa_1 k_t\le\dot{k}_t $. Choose $\kappa_2=\delta_H+b$, then based on assumption (\ref{ass1}), we conclude $-\kappa_2 h_t\le \dot{h}_t$.  Since $|s_t|,|i_t|\le 1$ and from assumption (\ref{ass2}), $|\beta(\cdot)|\le 1$, choose $\kappa_3=b+1$. Then $-\kappa_3 s_t\le \dot{s}_t$. From assumption (\ref{ass3}), $|\gamma(\cdot)|\le 1$, hence choose $\kappa_4=b+1$ and we have $-\kappa_4 i_t\le \dot{i}_t$.  Finally choose $\kappa_5=\delta_E$ then $-\kappa_5 e_t\le \dot{e}_t$. Now choose $F(y_t)=(-\kappa_1 k_t, -\kappa_2 h_t, -\kappa_3 s_t, -\kappa_4 i_t,  -\kappa_5 e_t)$. It is clear that $F(\cdot)$ is continuous on $\mathbb{R}_+^5$ and satisfy $ F(y_t)\le \dot{y}_t$
\end{proof}

\begin{lemma} \label{lemmaG}
Let $y_t=(k_t, h_t, s_t, i_t, e_t)$ and $x_t=(c_t, m_t, A_t)$ be state and control variables given in social planner problem (\ref{social}). Let $\bar{y_t}=(\bar{y}_1, \cdots,\bar{y}_5)$ an external variable  Then there exists a continuous function $G(\cdot):  \mathbb{R}_{+}^5\times\mathbb{R}_{+}^5\times \mathbb{R}_{+}^3 \to \mathbb{R}^5$ such that \[\dot{y}_t \le G(\bar{y}_t,y_t,x_t)\] and $G^j$ is concave with respect to $x_t$
\end{lemma}
\begin{proof}
Choose $\bar{y_t}=(1, \cdots,1)$ and $G(\cdot)=(f(k_t,1-i_t),g(m_t),b,1,E(A_t,e_t))$. We claim $G(\cdot)$ is continuous on $ \mathbb{R}_{+}^5\times\mathbb{R}_{+}^5\times \mathbb{R}_{+}^3$, satisfy $\dot{y}_t \le G(\bar{y}_t,y_t,x_t)$ and $G^j$ is concave with respect to $x_t$ for all $j\in\{\,1,2,\cdots,5\}\,$\\
\indent It is obvious that $G(\cdot)$ is continuous and satisfy $\dot{y}_t \le G(\bar{y}_t,y_t,x_t)$. Except $g(m_t)$ and $E(A_t,e_t)$ all other $G^j$ are constants respect to $x_t$, hence concave. For $g(m_t)$, under section \ref{assumptions} assumption (\ref{ass5}), concavity follows respect to $x_t$. Similarly concavity of $E(A_t,e_t)$ with respect to $x_t$ follows from assumption (\ref{ass8}).
\end{proof}

\begin{lemma}\label{lemmaY}
Let $y_t=(k_t, h_t, s_t, i_t, e_t)$ be the state variables given in social planner problem in (\ref{social}). There exists $a_j\ge 0, A_j>0$ such that \[\bar{y}_t^j, y_t^j\le A_je^{a_jt} \text{ and } \dot{y}^j_t\le A_je^{a_j t}\] for all $j\in \{\,1,2,\cdots, 5\}\,$
\end{lemma}
\begin{proof} First note that $\bar{y}_t^j=1$ for all $j$ from lemma (\ref{lemmaG}). Now, following proof of lemma 1 (with few simple modification) given in \cite{goenka2014infectious}, it can be easily shown that for any $a_1,a_2 \in (0,\theta)$ there exist $A_1, A_2>0$ such that $\max \{\,1, k_t, |\dot{k}_t|\}\,\le A_1e^{a_1 t}$ and $\max\{\,1, h_t, |\dot{h}_t|\}\,\le A_2e^{a_2t}$ for all $t$, where $\theta$ denote the discount rate used in social planner problem in (\ref{social}). Observe that $s_t,i_t,\beta\le 1$ and $|\dot{s}_t|\le |b|+|b|+1$.  Hence for any $a_3\in (0,\theta)$ there exist $A_3>0$ such that $\max\{\,1, s_t,|\dot{s}_t|\}\,\le A_3e^{a_3 t}$ for all $t$. Similar type arguments can be used to show that for any $a_4\in (0,\theta)$ there exist $A_4>0$ such that $\max\{\,1, i_t, |\dot{i}_t|\}\,\le A_4e^{a_4 t}$ for all $t$. 

Under assumption \ref{ass8}, we assumed $E_1(A_t,e_t)=\pdv{E(A_t,e_t)}{A_t}>0$ and $A_t\le A$ for all $t$. Hence $E(A_t,e_t)\le E(A,e_t)$. We also assumed that, $\lim_{e_t\to \infty}E_2(A_t,e_t)=\lim_{e_t\to \infty}\pdv{E(A_t,e_t)}{e_t}=0$. Hence from L'hopital's rule it follow that $\lim_{e_t\to \infty} \frac{E(A,e_t)}{e_t}=0$. Thus for any $a_5\in (0,\theta)$ there exists constants $C_0$ such that $\dot{e}_t=E(A_t,e_t)\le E(A,e_t)\le C_0+a_5e_t$. Again using similar argument given in first portion of lemma 1 in \cite{goenka2014infectious} , we conclude for any $a_5 \in (0, \theta)$ there exist $A_5>0$ such that $\max\{\,1,e_t, |\dot{e}_t|\}\,\le A_5e^{a_5 t}$. Hence the result.
\end{proof}

\begin{lemma} \label{lemmaX}
Let $x_t=(c_t, m_t, A_t)$ be control variables given in social planner problem \ref{social}. Then there exists $b_i>0, B_i>0$ for $i\in \{\,1,2,3\}\,$ such that $x_t^i\le B_ie^{b_i t}$
\end{lemma}
\begin{proof}
We know from assumption (\ref{ass4}), physical production function $f$ is increasing on $l_t=1-i_t$. Since $l_t\le 1$, we conclude that $f(k_t,l_t)\le f(k_t,1)$. Now using the fact that $\lim_{k_t\to \infty}f(k_t,1)=0$ and L'hopital's rule we conclude for any $b_0\in (0,\theta)$ there exists $C_0$ such that $f(k_t,l_t)\le f(k_t,1)\le C_0+b_0k_t$. But from lemma (\ref{lemmaY}), we know for any $b_0\in (0,\theta)$ there exist $A_0$ such that $k_t\le A_0e^{b_0 t}$. This imply  for any $b_0 \in (0, \theta)$ there exists $B_0$ such that $f(k_t,l_t)\le C_0+b_0k_t\le B_0e^{b_0 t}$. Now observe $m_t, c_t, A_t \le f(k_t,l_t)\le B_0e^{b_0 t}$ for all $t$. Hence the result.
\end{proof}
Before we continue observe that discount rate $\theta>b_i$ given in lemma \ref{lemmaX}. This is because, we choose $b_i$ such that $b_i \in (0, \theta)$.   Now we are ready to show the existence of a solution
\begin{theorem}
There exists a solution to the social planner's problem \ref{social}
\end{theorem}
\begin{proof}
Define $u(x_t)=u(c_t)$, where $u(c_t)$ satisfy the assumption \ref{ass6}. Hence $u(x_t)$ satisfy the { assumption } A5. From lemma \ref{lemmaF} to  \ref{lemmaX}, we can see that problem \ref{social} satisfy assumptions A1-A4. Hence from the theorem \ref{exsistancemain}, result follows.
\end{proof}

\section{Steady States and Social Planner Problem}
In this section we investigate the steady states of the social planner problem at optimality. We use first order necessary condition for optimal solution to find the steady states. For SIR model with the vaccination  given in Equation (\ref{sird}) , we found two steady states. One related to the disease free and the other related to the endemic state.  Economic steady states solutions can be divide based on those epidemic disease free and endemic steady states to have more clear understanding of how economy and epidemic interact with each other. Before we continue, we would suppress any variable in the form of $x_t$ into $x$ to make equations more readable. 

Before we write the current value Lagrangian for the social planner problem given in Equation (\ref{social}), observe that capital and labor both are essential to the physical production function hence $k>0, l>0$. If $m=0$ then using Equation (\ref{eqh}), we can conclude at the steady state, $h=0$ and if $m>0$ then $h>0$. Thus $m\ge 0$ would imply $h\ge 0$. For existing knowledge, if $A=0$ then $e=0$. This imply no direct investment to control the disease and no existing knowledge to control the disease. We attached this possibility with disease free steady state. Any other non trivial steady state solutions, we attached with endemic steady states.  The current value Lagrangian for this social planner problem is given by
\begin{align}
\begin{split}
    L(\cdot) &=u(c_t) \\ 
    &+\lambda_1\big( f(k,1-i)-m-A-c-(\delta_K+b-\mu)k\big) \\ 
    &+\lambda_2\big(g(m)-(\delta_H+b-\mu )h\big) \\
    &+\lambda_3\big((1-p)b-b s-\beta si\big)\\
    &+\lambda_4\big(\beta s i-(\gamma+b)i\big) \\
    &+\lambda_5(E-\delta_Ee)\\
    &+\nu_1i+\nu_2 m+ \nu_3 A +\nu_4 s
    \end{split}
\end{align}
where $\lambda_i$ are co-state variables and $\nu_i$ are Lagrangian multipliers. The first order necessary condition for optimal solution is given by the theorem 14.5 in \cite{caputo2005foundations}. First order conditions are given by,  

\small
\begin{align}
\pdv{L}{c} &=u'(c)-\lambda_1=0 \label{eqc} \\
\pdv{L}{m} &=-\lambda_1+ \lambda_2 g'(m)+\nu_2=0 \label{eqhm}\\ 
{\pdv{L}{A} } &=-\lambda_1-\lambda_3 \beta_1 s i +\lambda_4 (\beta_1 s i-\gamma_1 i)+\lambda_5 E_1+\nu_3=0 \label{eqA}\\
\dot{\lambda}_1 &=-\pdv{L}{k}+\lambda_1 \theta= \lambda_1(\theta+\delta_K+b-\mu-f_1)\label{eqlam1}\\
\dot{\lambda}_2 &=-\pdv{L}{h}+\lambda_2 \theta=\lambda_2 (\theta+b+\delta_H-\mu) +\lambda_3 \beta_3 si  -\lambda_4(\beta_3 si-\gamma_3i)\\
\dot{\lambda}_3 &=-\pdv{L}{s}+\lambda_3 \theta=\lambda_3(\theta+b+ \beta i)-\lambda_4 \beta i -\nu_4 \label{eqlam3}\\
\dot{\lambda}_4 &=-\pdv{L}{i}+\lambda_4 \theta= \lambda_1f_2+\lambda_3 \beta s+\lambda_4(\theta+b +\gamma-\beta s)-\nu_1 \label{eqlam4}\\
\dot{\lambda}_5 &=-\pdv{L}{e}+\lambda_5 \theta= \lambda_3 \beta_2 s i-\lambda_4(\beta_2 si-\gamma_2i)+\lambda_5(\theta+\delta_E-E_2) \label{eqE}\\
\pdv{L}{\lambda_1} &=\dot{k}_t, \pdv{L}{\lambda_2} =\dot{h}_t, \pdv{L}{\lambda_3} =\dot{s}_t, \pdv{L}{\lambda_4} =\dot{i}_t, \pdv{L}{\lambda_5} =\dot{e}_t\\
lim_{t \to\infty} &e^{-\theta t} \lambda_{1,t} k_t =0, lim_{t\to\infty}e^{-\theta t} \lambda_{2,t} h_t=0, lim_{t\to\infty}e^{-\theta t} \lambda_{3,t} s_t=0 \\
lim_{t\to\infty} &e^{-\theta t} \lambda_{4,t} i_t=0, lim_{t\to \infty}e^{-\theta t} \lambda_{5,t} e_t=0 \\
 \nu_i&\ge 0,  i\nu_1=0,  m\nu_2=0,  A\nu_3=0, s\nu_4=0 \label{eqconst}
\end{align}
\normalsize
 Let $x_t=(c_t,m_t, A_t),y_t=(k_t, h_t, s_t, i_t, e_t)$,  $\lambda_t=(\lambda_1, \lambda_2,\lambda_3,\lambda_4,\lambda_5)$ and $\nu_t=(\nu_1,\nu_2,\nu_3,\nu_4)$. Then steady state for social planner problem is given by a set of values $(x^*_t,y^*_t,\lambda^*_t,\nu^*_t)$ which  satisfy $\dot{y}_t=\dot{\lambda}_t=0$. Before we continue, we make following assumption regarding $E_2(A,e)$. 
 \begin{assumption}\label{assE2}
  For all $A, e$ we have $E_2(A,e)\ne \theta+\delta_E$
 \end{assumption}

 \begin{proposition}\label{propDF}
 There exist a disease free steady state with $m^*=0, h^*=0, s^*=1-p, i^*=0,  e^*=0, A^*=0$,  $k^*$ is determined by 
 \[ f_1(k^*,1)=(\theta+\delta_k +b-\mu)\] and $c^*$ is given by
 \[c^*=f(k^*,1)-(\delta_K+b-\mu)k^*\]
 \end{proposition}
\begin{proof}
Under the disease free steady state, $s^*=1-p$ and $i^*=0$. Hence $l^*=1-i^*=1$ and $\nu_4=0$.  From Equations (\ref{eqc}-\ref{eqhm}), we can see $\nu_2=u'(c)-\lambda_2g'(m)$. {When $l^*=1$, $\dot{\lambda}_2=\lambda_2(\theta+b+\delta_H-\mu)$}. Since $b\ge \mu$, we conclude that steady state value, $\lambda_2^*=0$. This imply $\nu_2=u'(c)>0$. Hence, $m^*=0$ (and $h^*=0$). From Equation (\ref{eqlam1}) we have $\lambda_1(\theta+\delta_K+b-\mu-f_1)=0$. However $\lambda_1^*=u'(c^*)>0$ from Equation (\ref{eqc}). Thus we conclude  there exist unique $k^*\ne 0$,  such that \[f_1(k^*,1)=(\theta+\delta_k +b-\mu)\] When $i^*=0$, Equation (\ref{eqE}) reduce to $\lambda_5(\theta+\delta_E-E^*_2)=0$. Since $E_2^*\ne \theta+\delta_E$ from assumption (\ref{assE2}), we conclude that $\lambda_5^*=0$. {When $i^*=0$ and using the fact $\nu_4=0$}, we can conclude from Equation (\ref{eqlam3}) that $\lambda_3^*=0$. From Equation (\ref{eqA}), we conclude that $\nu_3=\lambda_1>0$. Hence, $A^*=0$. Now using the fact $E(0,e)=0$ and the Equation (\ref{eqe}), $e^*=0$. Observe that $\nu_1$ is a non-negative free variable. Thus, we can solve for $\lambda_4$ using Equation (\ref{eqlam4}).  
\end{proof}

\begin{proposition}\label{prop43}
Steady states solution for endemic state exist if and only if there exists solution $m^*,k^*, A^*,e^*, h^*$ to following set of equations and inequalities
\begin{align}
(1-p) \beta(A,e,h) &>b+\gamma(A,e,h)\label{eqendemic}\\
     g(m)&=(b+\delta_H-\mu)h \label{eqgm}\\
     f_1(k,l) &=(\theta+\delta_K+b-\mu) \label{eqfkl}\\
     e &\ge 0 \label{eqepos}\\
     E(A,e) &=\delta_Ee \label{eqAe}\\
    m &\ge 0 \label{eqmnu21}\\
    (\theta+b+\delta_H-\mu) &\ge  l_{\theta,3}(A,e,h)f_2g'(m)\label{eqmnu22}\\
      0& =m\bigg( (\theta+b+\delta_H-\mu) -  l_{\theta,3}(A,e,h)f_2g'(m)\bigg)\label{eqnu23} \\
    A &\ge 0\label{eqAnu31}\\
 (\theta+\delta_E-E_2) &\ge l_{\theta,1}(A,e,h)f_2(\theta+\delta_E-E_2)+l_{\theta,2}(A,e,h)f_2 E_1 \label{eqAnu32}\\
 0 &= A\bigg((\theta+\delta_E-E_2) - l_{\theta,1}(A,e,h)f_2(\theta+\delta_E-E_2)-l_{\theta,2}(A,e,h)f_2 E_1\bigg) \label{eqAnu3}
\end{align}
Where
\[ l_{\theta,j}(A,e,h) =\frac{i\bigg(\gamma_j\beta i-(\beta_js-\gamma_j)(\theta+b)\bigg)}{\bigg((\theta+b+\gamma)(\theta+b+\beta i)-\beta s (\theta+b)\bigg)}\]
\end{proposition}
\begin{proof}

Let $m^*, h^*, s^*, i^*,  e^*, A^*$,  $k^*$ be a solution set to above system of equations and inequalities. We claim these values satisfy Equations (\ref{eqk}-\ref{eqineq}) and (\ref{eqc}-\ref{eqconst}) at steady state.  Since $(1-p) \beta(A^*,e^*,h^*) >b+\gamma(A^*,e^*,h^*)$ we know these $A^*,e^*,h^*$ corresponds to endemic state of epidemic model.  Denote $\beta(A^*,e^*, h^*)$ and $\gamma(A^*,e^*, h^*)$ by $\beta^*$ and $\gamma^*$. Then endemic steady state epidemic model values are given by, $s^*=\frac{\gamma^*+b}{\beta^*}, i^*=\frac{(1-p)b}{\gamma^*+b}-\frac{b}{\beta^*}$ and effective labor at endemic steady state is given by 
\begin{equation}
l^*=1-i^*=1-\bigg(\frac{(1-p)b}{\gamma^*+b}-\frac{b}{\beta^*}\bigg) \label{eqlabor}
\end{equation}
with $\dot{s}_t=\dot{i}_t=0$. Since $s^*, i^*>0$, we have $\nu_1=0, \nu_4=0$. Since $m^*$ and $h^*$ satisfy Equation (\ref{eqgm}), from Equation (\ref{eqh}) we conclude $\dot{h}_t=0$. Since, $k^*$ satisfy Equation (\ref{eqfkl}), we have 
\begin{equation}
c^*=f(k^*,l^*)-m^*-A^*-(\delta_K+b-\mu)k^*
\end{equation}
Thus, from equation ($\ref{eqk}$), we conclude that $\dot{k}_t=0$. 
Since $l^*<1$,  $k^*$ would be different from  $k^*$ in disease free steady state. Since $A^*,e^*$ satisfy Equation (\ref{eqepos}-\ref{eqAe}), we have $\dot{e}_t=0$. {Therefore, } if $m^*, h^*, s^*, i^*,  e^*, A^*$,  $k^*$ is a solution to above set of equations and inequalities, then $\dot{k}_t=\dot{h}_t=\dot{s}_t=\dot{i}_t=\dot{e}_t=0$. Next observe that since $m^*, h^*, s^*, i^*,  e^*, A^*$,  $k^*$ satisfy Equation (\ref{eqfkl}), from Equation (\ref{eqlam1}), we conclude that $\dot{\lambda}_1=0$

Define, 
\begin{alignat}{3}\label{eqm}
    m_{11} &=g'(m^*) \quad
    &m_{22} &=-\beta_1^* s^* i^*\quad
    &m_{23} &=\beta_1^* s^* i^*-\gamma_1^* i^*\\ \nonumber
    m_{24} &=E_1^*\quad
   & m_{31} &=\theta+b+\delta_H-\mu\quad
    & m_{32} &=\beta_3^* s^*i^*\\\nonumber
    m_{33} &=-(\beta_3^* s^*i^*-\gamma_3^*i^*) \quad
    &m_{42} &=\theta+b+ \beta^* i^*\quad
    &m_{43} &=-\beta^* i^*\\\nonumber
    m_{52} &=\beta^* s^*\quad
    &m_{53} &=\theta+b +\gamma^*-\beta^* s^*\quad
    &m_{62} &=\beta_2^* s^* i^*\\\nonumber
    m_{63} &=-(\beta_2^* s^*i^*-\gamma_2^*i^*) \quad
    &m_{64} &=(\theta+\delta_E-E_2^*)\quad
    &b_1 &=u'(c^*) \\ \nonumber
    b_2 &=u'(c^*) \quad
&b_5 &=-u'(c^*)f_2^*  \nonumber
\end{alignat}
Then we can write the system of first order conditions as follows
\begin{alignat*}{3}
A &=
    \begin{bmatrix}
    m_{11} & 0 & 0 & 0 & 1 & 0 \\
    0 & m_{22} & m_{23} & m_{24}& 0 & 1\\
     m_{31} & m_{32} & m_{33} & 0 & 0 & 0\\
     0& m_{42} & m_{43} & 0 &0 &0 \\
     0& m_{52}& m_{53} &0 &0 &0 \\
     0 & m_{62} & m_{63} & m_{64} & 0 &0 
    \end{bmatrix}
    &B&=\begin{bmatrix}
    \lambda_2 \\ \lambda_3\\ \lambda_4\\ \lambda_5\\ \nu_2\\ \nu_3
    \end{bmatrix} 
    &C&=
    \begin{bmatrix}
    b_1 \\ b_2 \\ 0\\ 0 \\b_5 \\ 0
    \end{bmatrix}
\end{alignat*}
The determinant of the matrix $A$ is given by \[|A|=m_{31}m_{64}(m_{42}m_{53}-m_{43}m_{52})\]
Observe that by assumption (\ref{eqE}), $m_{64}\ne 0$. Since $b\ge \mu$, we  have $m_{31}>0$ . Hence, if $(m_{42}m_{53}-m_{43}m_{52})$ is nonzero, then $A$ is invertible.
Note that 
\begin{align*}
   m_{42}m_{53}-m_{43}m_{52} &= (\theta+ b+\beta i)(\theta +b+\gamma-\beta s)-(-\beta i)(\beta s) \\
    &=  (\theta+b+\gamma)(\theta+b+\beta i)-\beta s (\theta+b)
\end{align*}
Since we assume Equation (\ref{eqendemic}) holds, we know endemic steady state exists and  $s=\frac{\gamma^*+b}{\beta^*}$. Hence,
$m_{53}m_{64}- m_{54}m_{63}>0$. Thus, endemic steady state values of $\lambda_i^*$'s and $\nu_2^*, \nu_3^*$ exists and given by $B=A^{-1}C$ and $\lambda_1^*=u'(c^*)$. Using  notations in \cite{goenka2014infectious}, define
 
\[l_{\theta,j}(A,e,h)=\frac{i\bigg(\gamma_j\beta i-(\beta_js-\gamma_j)(\theta+b)\bigg)}{\bigg((\theta+b+\gamma)(\theta+b+\beta i)-\beta s (\theta+b)\bigg)}\]
where $j=1,2,3$. Then solutions for $\nu_2$ and $\nu_3$ are given by
\begin{align}
 \nu_2 & = b_1+\frac{b_5\bigg(m_{42}m_{33}-m_{43}m_{32}\bigg)m_{11}}{m_{31}\bigg(m_{42}m_{53}-m_{43}m_{52}\bigg)}\\
 &=u'(c)-u'(c)\frac{l_{\theta,3}(A,e,h)}{\big(\theta+b+\delta_H-\mu\big)}f_2g'(m)\label{eqnu2}\\
 \nu_3 &= b_2 +b_5\frac{\bigg(m_{22}m_{43}-m_{23}m_{42}\bigg)}{\bigg(m_{42}m_{53}-m_{43}m_{52}\bigg)}+b_5\frac{m_{24}\bigg(m_{63}m_{42}-m_{43}m_{62}\bigg)}{m_{64}\bigg(m_{42}m_{53}-m_{43}m_{52}\bigg)}\\
 &= u'(c) -u'(c)l_{\theta,1}(A,e,h)f_2-u'(c)l_{\theta,2}(A,e,h)f_2\bigg(\frac{E_1}{\theta+\delta_E-E_2}\bigg)\label{eqnu3}
\end{align}
Since $A^*,e^*,h^*$ satisfy Equations (\ref{eqmnu21}-\ref{eqAnu3}), we conclude that Equations (\ref{eqconst}) holds. Now observe that for these $\lambda_i$ values we have $\pdv{L}{c}=\pdv{L}{m}=\pdv{L}{A}=\dot{\lambda}_i=0$. Hence this solution satisfy first order conditions too.  \\

Conversely suppose $m^*, h^*, s^*, i^*,  e^*, A^*$,  $k^*$  is solution to endemic steady state for equations given in (\ref{eqk}-\ref{eqineq}) and satisfy first order conditions given in (\ref{eqc}-\ref{eqconst}). Then, clearly we can see that it need to satisfy the equations and inequalities given in (\ref{eqc}-\ref{eqconst}) 
\end{proof}
We make following assumption in order to further characterize endemic steady state.

\begin{assumption}\label{assl} 
{Assume } $l_{\theta,j}(A,e,h)$ decrease in $A,e$ and $h$. That is $\pdv{l_{\theta,j}(A,e,h)}{A}<0,\pdv{l_{\theta,j}(A,e,h)}{e}<0, \pdv{l_{\theta,j}(A,e,h)}{h}<0 $ for $j=1,2,3$
\end{assumption}

\begin{proposition}\label{lthetaprop}
$l_{\theta,j}(0,0,0)$ decrease with respect to $\theta$ for all $j=1,2,3$
\end{proposition}
\begin{proof}
First observe that, 
\begin{align}
    l_{\theta,j}(0,0,0) &=\frac{\bar{i}\bigg(\gamma_j(0,0,0)\bar{\beta} \bar{i}-(\beta_j(0,0,0)\underline{s}-\gamma_j(0,0,0))(\theta+b)\bigg)}{\bigg((\theta+b+\underline{\gamma})(\theta+b+\bar{\beta}\bar{i})-\bar{\beta} \underline{s} (\theta+b)\bigg)} \nonumber
    \end{align}
    Hence,
    \begin{align}
   \bigg(\frac{1}{-\bar{i}(\beta_j\underline{s}-\gamma_j)}\bigg)\dv{l_{\theta,j}(0,0,0)}{\theta} &=\frac{-\theta^2-2(b+\eta)\theta+(\underline{\gamma}-\eta)\bar{\beta}\bar{i}-b(\eta+b)}{\bigg((\theta+b+\underline{\gamma})(\theta+b+\bar{\beta}\bar{i})-\bar{\beta} \underline{s} (\theta+b)\bigg)^2} \label{theta1}
\end{align}
where $\underline{s}=\frac{\underline{\gamma}+b}{\bar{\beta}}, \bar{i}=\frac{(1-p)b}{\underline{\gamma}+b}-\frac{b}{\bar{\beta}}$ and $\eta=\frac{\gamma_j\bar{\beta}\bar{i}}{-(\beta_j\underline{s}-\gamma_j)}$\\
Observe that numerator of Equation (\ref{theta1}) is a quadratic polynomial in $\theta \in (0,\infty)$. It is easy to show that constant term of that polynomial is negative. hence the result follows.
\end{proof}

In next few propositions, we further {analyze} the endemic steady state. We denote effective labor when $A=0,e=0,h=0$ as $\underline{l}=1-\big(\frac{(1-p)b}{\underline{\gamma}}-\frac{b}{\bar{\beta}}\big)$. For given $b, \theta $ define $\underline{k}, \bar{k}$ such that   $f_1(\underline{k},\underline{l})=(\theta+\delta_K +b-\mu)$, $f_1(\bar{k},1)=(\theta+\delta_K +b-\mu)\bar{k}$. For each fixed $b\in[\mu,(1-p)\bar{\beta}-\underline{\gamma})$, let $\hat{\theta}_1(b),\hat{\theta}_2(b)$  to be solutions to following equations: 
\begin{alignat}{1}
\hat{\theta}_1(b)&: (\theta+b+\delta_H-\mu) =  l_{\theta,3}(0,0,0)f_2(\underline{k},\underline{l})g'(0) \label{t1}\\
\hat{\theta}_2(b)&:  (\theta+\delta_E-\bar{E}_2) = l_{\theta,1}(0,0,0)f_2(\underline{k},\underline{l})(\theta+\delta_E-\bar{E}_2)+l_{\theta,2}(0,0,0)f_2(\underline{k},\underline{l}) \bar{E}_1 \label{t2}
\end{alignat}
\begin{proposition}
For each $b\in[\mu,(1-p)\bar{\beta}-\underline{\gamma})$, there exist unique solution $\hat{\theta}_1(b),\hat{\theta}_2(b)$ to above Equations (\ref{t1}-\ref{t2})
\end{proposition}
\begin{proof}
For Equation (\ref{t1}) it is clear that LHS is increasing w.r.t to $\theta$ and RHS is decreasing w.r.t to $\theta$ without any asymptotes. Hence, there exist unique solution $\hat{\theta}_1(b)$.  Rewrite Equation (\ref{t2}) as below. \[1 = l_{\theta,1}(0,0,0)f_2(\underline{k},\underline{l})+l_{\theta,2}(0,0,0)f_2(\underline{k},\underline{l}) \frac{\bar{E}_1}{(\theta+\delta_E-\bar{E}_2)} \] Then observe that RHS  of above equation is decreasing w.r.t $\theta>\bar{E}_2-\delta_E$ and LHS is  a constant w.r.t $\theta$. Also note that $\theta=\bar{E}_2-\delta_E$ is a vertical asymptote. Therefore, there exist a unique solution $\hat{\theta}_2(b)$. Note that, it is possible these solutions, $\hat{\theta}_1(b),\hat{\theta}_2(b)$ to be negative.
\end{proof}

We define $\hat{\theta}_{max}(b)=\max(\hat{\theta}_1(b),\hat{\theta}_2(b),0)$ and $\hat{\theta}_{min}(b)=\min(\hat{\theta}_1(b),\hat{\theta}_2(b))$ . 

\begin{proposition}\label{propbehave}
Let $m^*,k^*,l^*$ be endemic steady state solutions. As $A,e,h$ increases 
\begin{enumerate}
    \item $g'(m^*(h))$ decreases
    \item $l^*(A,e,h)$ increases
    \item $k^*(A,e,h)$ increases 
    \item  $f_2(k^*,l^*)$ decreases 
    \item $f_1(k^*,l^*)$ remain constant
\end{enumerate}
\end{proposition}
\begin{proof}
From Equations (\ref{eqgm},\ref{eqfkl},\ref{eqlabor}) it is clear that solution $m^*,k^*,l^*$ can be written as functions of $A,e,h$. First we claim that that $g'(m^*(h))$ decrease as $h$ increases. Observe that from Equation (\ref{eqgm}), 
\begin{align}
    \pdv{g'(m^*(h))}{h} &=g''(m^*)\pdv{m^*(h)}{h} \nonumber\\
    &=g''(m^*) \bigg(\frac{b+\delta_H-\mu}{g'(m^*)}\bigg)
\end{align}
Since $g'(m^*)>0, g''(m^*)<0$ from assumption (\ref{ass5}), we conclude that $\pdv{g'(m^*)}{h}<0$. Hence $g'(m^*(h))$ decrease as $h$ increases. Next, observe that at endemic steady state, 
\begin{align}
\pdv{l^*(A,e,h)}{A} &=\frac{(1-p)b\gamma_1}{(\gamma +b)^2}-\frac{b\beta_1}{\beta^2}
\end{align}
Since $\gamma_1>0$ and $\beta_1<0$ from assumptions (\ref{ass2}-\ref{ass3}). we conclude that $\pdv{l^*(A,e,h)}{A}>0$. Similarly $\pdv{l^*(A,e,h)}{e}>0$ and $\pdv{l^*(A,e,h)}{h}>0$. Hence $l^*(A,e,h)$ increases as $A,e,h$ increases.  From Equation (\ref{eqfkl}), \begin{align}
    \pdv{f_1}{k}\pdv{k^*}{A}+\pdv{f_1}{l}\pdv{l^*}{A} &=0\nonumber\\
    f_{11}\pdv{k^*}{A}+f_{12}\pdv{l^*}{A} &=0\nonumber\\
    \pdv{k^*(A,e,h)}{A} &=-\frac{f_{12}}{f_{11}}\pdv{l^*}{A}\nonumber
\end{align}
Since $f_{12}>0,f_{11}<0, \pdv{l^*}{A}>0$, we have $ \pdv{k^*}{A}>0$. Similar proof apply w.r.t $e$ and $h$. Thus, $k^*(A,e,h)$ increases as $A,e,h$ increases. Now, we want to see the behaviour of $f_2(k^*,l^*)$ as $A,e,h$ increases. 
\begin{align}
    \pdv{f_2(k^*,l^*)}{A} &=\pdv{f_2}{k}\pdv{k^*}{A}+\pdv{f_2}{l}\pdv{l^*}{A} \nonumber \\
    &=f_{21}\bigg(-\frac{f_{12}}{f_{11}}\bigg)\pdv{l^*}{A}+f_{22}\pdv{l^*}{A}\nonumber \\
    &=\bigg(\frac{f_{22}f_{11}-f_{12}f_{21}}{f_{11}}\bigg)\pdv{l^*}{A}\nonumber 
\end{align}
Since $f_{11}<0,f_{22}f_{11}-f_{12}f_{21}>0, \pdv{l^*}{A}>0$, we conclude $f_{2}(k^*,l^*)$ decreases w.r.t $A$. Similarly, $f_2(k^*,l^*)$ decreases w.r.t $e,h$

Next,from Equation (\ref{eqfkl}),
\begin{align}
    \pdv{f_1(k^*,l^*)}{A} &=0\nonumber 
\end{align}
Thus, $f_{1}(k^*,l^*)$ remain constant w.r.t $A$. Similarly, $f_1(k^*,l^*)$ remains constant w.r.t $e,h$
\end{proof}

\begin{proposition} \label{propEndemicNo}
Under assumptions (\ref{ass1}-\ref{assl}) , for each fixed $b\in[\mu,(1-p)\bar{\beta}-\underline{\gamma})$ if $\theta > \hat{\theta}_{max}(b)$, then  $A^*=0, e^*=0, ^*h=0, m^*=0, k^*=\underline{k},l^*=\underline{l},c^*= f(\underline{k},\underline{l})-\big(\delta_K+b-\mu\big)\underline{k}$ is the unique steady state solution to the system given in Equations (\ref{eqendemic}-\ref{eqAnu3}). In otherwords, when $\theta > \hat{\theta}_{max}$ there exist unique endemic steady state without health expenditure and without direct investment to control the epidemic.
\end{proposition}

\begin{proof}
Since $b\in [\mu,(1-p)\bar{\beta}-\underline{\gamma})$ we know, we are in endemic state. Let $A^*=0, e^*=0,h^*=0, m^*=0$, $ k^*=\underline{k},l^*=\underline{l}$ and $\theta>\hat{\theta}_{max}(b)$. We want to claim that this satisfy equations and inequalities in (\ref{eqendemic}-\ref{eqAnu3}). Note that only conditions we need to check are inequalities (\ref{eqmnu22}) and (\ref{eqAnu32}). Observe that, $(\theta+b+\delta_H-\mu)> (\hat{\theta}_{max}(b)+b+\delta_H-\mu)\ge (\hat{\theta}_{1}(b)+b+\delta_H-\mu)$ and $l_{\hat{\theta}_1,3}(0,0,0)f_2(\underline{k},\underline{l})g'(0)\ge l_{\hat{\theta}_{max},3}(0,0,0)f_2(\underline{k},\underline{l})g'(0)>l_{\theta,3}(0,0,0)f_2(\underline{k},\underline{l})g'(0)$. Hence inequality (\ref{eqmnu22}) is satisfied. Also note that, 
\begin{align*}
1 &= l_{\hat{\theta}_2,1}(0,0,0)f_2(\underline{k},\underline{l})+l_{\hat{\theta}_2,2}(0,0,0)f_2(\underline{k},\underline{l}) \frac{\bar{E}_1}{(\hat{\theta}_2(b)+\delta_E-\bar{E}_2)}\\
&>l_{\hat{\theta}_{max},1}(0,0,0)f_2(\underline{k},\underline{l})+l_{\hat{\theta}_{max},2}(0,0,0)f_2(\underline{k},\underline{l}) \frac{\bar{E}_1}{(\hat{\theta}_{max}(b)+\delta_E-\bar{E}_2)}\\
&>l_{\theta,1}(0,0,0)f_2(\underline{k},\underline{l})+l_{\theta,2}(0,0,0)f_2(\underline{k},\underline{l}) \frac{\bar{E}_1}{(\theta+\delta_E-\bar{E}_2)}
\end{align*}
Hence inequality (\ref{eqAnu32}) is satisfied. Hence $A^*=0, e^*=0,h^*=0, m^*=0$, $ k^*=\underline{k},l^*=\underline{l}$ is a solution to endemic steady state when $\theta>\hat{\theta}_{max}(b)$. The solution become unique from proposition (\ref{propbehave}). Let $\theta>\hat{\theta}_{max}$ and $A^*,e^*,h^*$ be a solution to the system given  (\ref{eqendemic}-\ref{eqAnu3}). Then,

$l_{\theta,3}(0,0,0)f_2(\underline{k},\underline{l})g'(0)>l_{\theta,3}(A^*,e^*,h^*)f_2(k^*,l^*)g'(m^*)$ and 
\begin{align*}
l_{\theta,1}(0,0,0)f_2(\underline{k},\underline{l})+l_{\theta,2}(0,0,0)f_2(\underline{k},\underline{l}) \frac{\bar{E}_1}{(\theta+\delta_E-\bar{E}_2)}\\
>l_{\theta,1}(A^*,e^*,h^*)f_2(k^*,l^*)+l_{\theta,2}(A^*,e^*,h^*)f_2(k^*,l^*) \frac{E_1(A^*,e^*)}{(\theta+\delta_E-E_2(A^*,e^*))}
\end{align*}

Thus, we conclude $\nu_2>0$ and $\nu_3>0$. Hence $A^*=0, m^*=0$, which imply, $e^*=0,h^*=0$. This would also imply $l^*=\underline{l}$ and $k^*=\underline{k}$. Hence solution is unique.
\end{proof}

\begin{proposition}\label{propEndemicPos_A}
Under assumptions (\ref{ass1}-\ref{assl}), for each fixed $b\in[\mu,(1-p)\bar{\beta}-\underline{\gamma})$ if $\hat{\theta}_{max}(b)=\hat{\theta}_2(b)>0$, then for all $\theta \in (\max\{\,0,\hat{\theta}_{1}(b)\}\,,\hat{\theta}_{2}(b))$ there exists a unique endemic steady state with zero health expenditure and nonzero epidemic control investment with solution determined by the following system of equations
\begin{align}
l &=1-\bigg(\frac{(1-p)b}{\gamma+b}-\frac{b}{\beta}\bigg)\\
       m &=0, h=0\\
     f_1(k,l) &=(\theta+\delta_K+b-\mu) \\
     E(A,e) &=\delta_Ee \\  
     (\theta+\delta_E-{E}_2(A,e))  &= l_{\theta,1}(A,e,0)f_2(k,l)(\theta+\delta_E-{E}_2(A,e)) \nonumber\\
     &+l_{\theta,2}(A,e,0)f_2(k,l) {E}_1(A,e) \label{eqAe3}\\
     c &=f(k,l)-A-(\delta_K+b-\mu)k
\end{align}
\end{proposition}
\begin{proof}
 Since $(1-p) \beta(A^*,e^*,h^*) >b+\gamma(A^*,e^*,h^*)$ for any $A^*,e^*$ and $h^*$, we know endemic steady state for epidemic model exist and $l^*$ is given by Equation (\ref{eqlabor}). Let $m^*,h^*,e^*,A^*,k^*$ be a solution to system of equations and inequalities given in (\ref{eqendemic}-\ref{eqAnu3}) . Since   $\hat{\theta}_{1}(b)<\theta<\hat{\theta}_{2}(b)$, 
\begin{align}
         (\theta+b+\delta_H-\mu)  &>  l_{\theta,3}(A^*,e^*,h^*)f_2(k^*,l^*)g'(m^*) 
\end{align}
 hence $m^*=0,h^*=0$. Suppose $A^*=0$ then $e^*=0$ and  $l^*=\underline{l}$ and $k^*=\underline{k}$. This would imply 
 \begin{align}
          (\theta_2(b)+\delta_E-\bar{E}_2) &> (\theta+\delta_E-\bar{E}_2)\\ l_{\theta,1}(0,0,0)f_2(\underline{k},\underline{l})(\theta+\delta_E-\bar{E}_2) &+l_{\theta,2}(0,0,0)f_2(\underline{k},\underline{l}) \bar{E}_1 \nonumber\\
          >l_{\theta_2(b),1}(0,0,0)f_2(\underline{k},\underline{l})(\theta_2(b)+\delta_E-\bar{E}_2) &+l_{\theta_2(b),2}(0,0,0)f_2(\underline{k},\underline{l}) \bar{E}_1
\end{align}
 We also know that,
  \begin{align}
           (\theta+\delta_E-\bar{E}_2) &\ge l_{\theta,1}(0,0,0)f_2(\underline{k},\underline{l})(\theta+\delta_E-\bar{E}_2) +l_{\theta,2}(0,0,0)f_2(\underline{k},\underline{l}) \bar{E}_1
\end{align}
 and 
   \begin{align}
           (\theta_2(b)+\delta_E-\bar{E}_2) &= l_{\theta_2(b),1}(0,0,0)f_2(\underline{k},\underline{l})(\theta_2(b)+\delta_E-\bar{E}_2) +l_{\theta_2(b),2}(0,0,0)f_2(\underline{k},\underline{l}) \bar{E}_1
\end{align}
 Hence only possibility is 
   \begin{align}
           (\theta+\delta_E-\bar{E}_2) &= l_{\theta,1}(0,0,0)f_2(\underline{k},\underline{l})(\theta+\delta_E-\bar{E}_2) +l_{\theta,2}(0,0,0)f_2(\underline{k},\underline{l}) \bar{E}_1
\end{align}
However, this is a contradiction to the definition of $\hat{\theta}_{2}(b)$ given in Equations (\ref{t2}). Therefore $A^*>0$  and $A^*$ is determined by the Equation (\ref{eqAe3}). 
{The  uniqueness of the solution follows from implicit mapping Theorem.} 
\end{proof}

\begin{proposition}
Under assumptions (\ref{ass1}-\ref{assl}), for each fixed $b\in[\mu,(1-p)\bar{\beta}-\underline{\gamma})$ if $\hat{\theta}_{max}(b)=\hat{\theta}_1(b)>0$, then for all $\theta \in (\max\{\,0,\hat{\theta}_{2}(b)\}\,,\hat{\theta}_{1}(b))$ there exists a endemic steady state with nonzero health expenditure and zero epidemic control investment with solution determined by the following system of equations
\begin{align}
l &=1-\bigg(\frac{(1-p)b}{\gamma+b}-\frac{b}{\beta}\bigg)\\
       A &=0, e=0\\
     f_1(k,l) &=(\theta+\delta_K+b-\mu) \\
     g(m)&=(b+\delta_H-\mu)h\\
     (\theta+b+\delta_H-\mu)  &=  l_{\theta,3}(0,0,h)f_2(k,l)g'(m)\\
     c &=f(k,l)-m-(\delta_K+b-\mu)k
\end{align}
\end{proposition}
\begin{proof}
omitted
\end{proof}



\newpage
\section{Sufficient conditions}
We have found solutions for steady states and we know these solutions satisfy first order necessary conditions for optimality given in Equations (\ref{eqc}-\ref{eqconst}). Now we claim these solutions are indeed optimal solutions. Observe that under Mangasarian sufficient conditions, Lagrangian function must be concave with respect to both state and control variable for all $t \in [0,\infty)$ over an open convex set containing all the admissible values of the state and control variables \cite{caputo2005foundations}. 
Define  $x_t=(k_t, h_t, s_t, i_t, e_t)$ as state variables, $z_t=(c_t, m_t, A_t)$ as control variable set and $\lambda_t=(\lambda_{1,t},\cdots,\lambda_{5,t})$. Then Hamiltonian is given by
\begin{align}
\begin{split}
    H(x_t,z_t,\lambda_t) &=u(c_t) \\ 
    &+\lambda_1\big( f(k_t,l_t)-m_t-A_t-c_t-(\delta_K+b-\mu)k_t\big) \\ 
    &+\lambda_2\big(g(m_t)-(\delta_H+b-\mu )h_t\big) \\
        &+\lambda_3\big((1-p)b-b s_t-\beta s_ti_t\big)\\
    &+\lambda_4\big(\beta s_t i_t-(\gamma+b)i_t\big) \\
    &+\lambda_5(E(A_t,e_t)-\delta_E e_t)
    \end{split}
\end{align}

with first order and transversality conditions given as,

\small
\begin{align}
\pdv{H}{c} &=u'(c)-\lambda_1=0 \label{foncc}\\
\pdv{H}{m} &=-\lambda_1+ \lambda_2 g'(m)=0 \\ 
{\pdv{H}{A} } &=-\lambda_1-\lambda_3 \beta_1 s i +\lambda_4 (\beta_1 s i-\gamma_1 i)+\lambda_5 E_1=0 \label{foncA}\\
\dot{\lambda}_1 &= \lambda_1(\theta+\delta_K+b-\mu-f_1) \label{fonck}\\
\dot{\lambda}_2 &=\lambda_2 (\theta+b+\delta_H-\mu) +\lambda_3 \beta_3 si  -\lambda_4(\beta_3 si-\gamma_3i)\\
\dot{\lambda}_3 &=\lambda_3(\theta+b+ \beta i)-\lambda_4 \beta i \\
\dot{\lambda}_4 &= \lambda_1f_2+\lambda_3 \beta s+\lambda_4(\theta+b +\gamma-\beta s)\\
\dot{\lambda}_5 &= \lambda_3 \beta_2 s i-\lambda_4(\beta_2 si-\gamma_2i)+\lambda_5(\theta+\delta_E-E_2) \\
\pdv{H}{\lambda_1} &=\dot{k}_t, \pdv{H}{\lambda_2} =\dot{h}_t, \pdv{H}{\lambda_3} =\dot{s}_t, \pdv{H}{\lambda_4} =\dot{i}_t, \pdv{H}{\lambda_5} =\dot{e}_t \label{focend}\\
lim_{t \to\infty} &e^{-\theta t} \lambda_{1,t} k_t =0, lim_{t\to\infty}e^{-\theta t} \lambda_{2,t} h_t=0, lim_{t\to\infty}e^{-\theta t} \lambda_{3,t} s_t=0 \label{foctrans1}\\
lim_{t\to\infty} &e^{-\theta t} \lambda_{4,t} i_t=0, lim_{t\to \infty}e^{-\theta t} \lambda_{5,t} e_t=0 \label{foctrans2}
\end{align}
\normalsize

Then from theorem 21.5 of \cite{simonblume}, the Hamiltonian is concave w.r.t to state and control variable if and only if the Hessian w.r.t to those variables is negative semi definite. Then from theorem 16.2 of \cite{simonblume}, a matrix is negative semi definite if and only if all principal minors are alternate in sign so that odd order ones are non-positive and even order ones are non-negative. It is unclear whether one of  order one principal minor,
\begin{align}
    H_{AA}(x_t,z_t,\lambda_t) &= -\lambda_3 (\beta_{11}s_ti_t)+\lambda_4(\beta_{11}s_ti_t-\gamma_{11}i_t)+\lambda_5(E_{11})
\end{align} 
is non negative or non-positive. Therefore, Mangasarian condition is not easy to apply. Arrow sufficiency conditions require, maximized Hamiltonian function to be a concave function of state variables for all $t\in[0,\infty)$ over an open convex set containing all admissible values of state variables.  Define the maximized Hamiltonian $M(\cdot)$ as in \cite{caputo2005foundations},
\begin{equation}
    M(t,x_t,\lambda_t)=\max_{z_t} H(t,x_t,z_t,\lambda_t)
\end{equation}

Consider disease free endemic steady state.  Under the disease free steady state, current value Hamiltonian will be reduced to 
\begin{align}
\begin{split}
    H(x_t,z_t,\lambda_t) &=u(c_t) \\ 
    &+\lambda_1\big( f(k,1)-m-A-c-(\delta_K+b-\mu)k\big) \\ 
    &+\lambda_2\big(g(m)-(\delta_H+b-\mu )h\big) \\
    &+\lambda_5(E-\delta_Ee)
    \end{split}
\end{align}

Observe that physical production function $f(k,l)$, health production function $g(m)$ and knowledge production $E(A,e)$ all follows neoclassical growth theory \cite{caputo2005foundations}. Hence, the steady state solution given in in the proposition (\ref{propDF}) is optimal (see Chapter 15 in \cite{caputo2005foundations}). Now consider the endemic steady state without health expenditure and without direct investment to control the epidemic. Then, current value Hamiltonian is given by 
\begin{align}
    H(x_t,z_t,\lambda_t) &=u(c_t) +\lambda_1\big( f(k,1-i)-c-(\delta_K+b-\mu)k\big) 
\end{align}
Thus, this is correspond to neoclassical steady states and hence the solution given in the proposition (\ref{propEndemicNo}) is optimal in a local neighborhood where $m=0, A=0,h=0,e=0$. Now consider the endemic state with either positive health expenditure or positive direct investment to control epidemic. We follows the proof given in \cite{goenka2014infectious} with few modifications.

\begin{assumption}\label{assboundk}
For all $l\in (0,1)$ with $f_1(k,l)>b-d+\delta$ there exist a maximum sustainable capital stock, $\hat{k}$, for the given $l$. This maximum sustainable capital stock can be obtained by solving $f(\hat{k},l)=(b-d+\delta)\hat{k}$
\end{assumption}

\begin{proposition}
Given endemic state, let $x_t^*=(k_t^*,h_t^*,s_t^*,i_t^*,e_t^*), z_t^*=(c_t^*,m_t^*,A_t^*), $ and $\lambda_t=(\lambda_{1,t},\cdots,\lambda_{5,t})$ denote a path with positive health expenditure and positive direct investment to control the epidemic which satisfies necessary condition given in Equations (\ref{foncc}-\ref{focend}) along with the transversality conditions given in (\ref{foctrans1}-\ref{foctrans2}) and initial condition $x_0=(k_0,h_0,s_0,i_0,e_0)$. Then this path is a locally optimal for the social planner problem given in (\ref{social})
\end{proposition}
 
 \begin{proof}
Now suppose there exist another path $x_t,z_t,\lambda_t$ with same initial condition $x_0$. We want to claim 
\[\int_0^{\infty} e^{-\theta t}u(c_t^*)dt \ge \int_{0}^{\infty} e^{-\theta t} u(c_t)dt\]
Notice that,
\begin{align}\label{int1}
&\int_{0}^{\infty} e^{-\theta t} [ H(x_t^*,z_t^*,\lambda_t)-H(x_t,z_t,\lambda_t)+<\dot{\lambda}_t-\theta \lambda_t,x_t^*-x_t>]dt  \nonumber\\ &= \int_{0}^{\infty} e^{-\theta t} [H(x_t^*,z_t^*,\lambda_t)-H(x_t^*,z_t,\lambda_t)]dt +\int_{0}^{\infty} e^{-\theta t} [H(x_t^*,z_t,\lambda_t)-H(x_t,z_t,\lambda_t)]dt \nonumber\\  &+ \int_0^{\infty} e^{-\theta t}<\dot{\lambda}_t-\theta \lambda_t,x_t^*-x_t>dt
\end{align}

\noindent Consider the first term of Equation (\ref{int1}), denote $\beta^*=\beta(A_t^*,e_t,h_t), \gamma^*=\gamma(A_t^*,e_t,h_t)$ and $ E^*=E(A^*_t,e_t)$ then,
\begin{align*}
  H(x_t^*,z_t^*,\lambda_t)-H(x_t^*,z_t,\lambda_t) &=  \big(u(c_t^*)-u(c_t)\big)-\lambda_1(c_t^*-c_t)-\lambda_1(m_t^*-m_t)-\lambda_1(A_t^*-A_t) \\ &+\lambda_2(g(m_t^*)-g(m_t)) -\lambda_3(\beta^*-\beta)s_t^*i_t^*+\lambda_4\big((\beta^*-\beta)s_t^*i_t^*-(\gamma^*-\gamma)i_t^* \big) \\
  &+\lambda_5(E^*-E)
\end{align*}

\noindent from concavity of $u(\cdot)$ and $g(\cdot)$,
$\frac{u(c_t)-u(c_t^*)}{c_t-c_t^*}\le u'(c_t)$ and $\frac{g(m_t)-g(m_t^*)}{m_t-m_t^*}\le g'(m_t)$. Also, since $x_t^*,y_t^*,\lambda_t$ satisfies first order conditions, $\lambda_{1,t}=u'(c_t^*)=\lambda_{2,t}g'(m_t^*)$. Hence,
\begin{align*}
  H(x_t^*,z_t^*,\lambda_t)-H(x_t^*,z_t,\lambda_t) &\ge  u'(c_t^*) (c_t^*-c_t)-\lambda_1(c_t^*-c_t)-\lambda_1(m_t^*-m_t)-\lambda_1(A_t^*-A_t) \\ &+\lambda_2 g'(m_t^*)(m_t^*-m_t) -\lambda_3(\beta^*-\beta)s_t^*i_t^*\\ &+\lambda_4\big((\beta^*-\beta)s_t^*i_t^*-(\gamma^*-\gamma)i_t^* \big)+\lambda_5(E^*-E)\\
  &\ge -\lambda_1(A_t^*-A_t)-\lambda_3(\beta^*-\beta)s_t^*i_t^* +\lambda_4\big((\beta^*-\beta)s_t^*i_t^*-(\gamma^*-\gamma)i_t^* \big) \\ &+\lambda_5(E^*-E)
\end{align*}
Now using concavity of $\beta,\gamma$ and $E$, 
\begin{align*}
  H(x_t^*,z_t^*,\lambda_t)-H(x_t^*,z_t,\lambda_t) &\ge  -\lambda_1(A_t^*-A_t)-(\lambda_3-\lambda_4)\beta_1^*s_t^*i_t^*(A_t^*-A_t)-\lambda_4\gamma_1^*i_t^*(A_t^*-A_t) \\ &+\lambda_5 E_1^*(A_t^*-A_t)\\
  &\ge 0
\end{align*}
for last inequality we used first order conditions.  Therefore the first term of the Equation (\ref{int1}),
\begin{equation}\label{int2}
    \int_{0}^{\infty} e^{-\theta t} [H(x_t^*,z_t^*,\lambda_t)-H(x_t^*,z_t,\lambda_t)]dt \ge 0
\end{equation}

Now consider the second term of the Equation (\ref{int1}), 
\begin{align*}
    H(x_t^*,z_t,\lambda_t)-H(x_t,z_t,\lambda_t) &= \lambda_{1,t} \dv{(k_t^*-k_t)}{t}+\lambda_{2,t} \dv{(h_t^*-h_t)}{t}+\lambda_{3,t} \dv{(s_t^*-s_t)}{t} \\ &+\lambda_{4,t}  \dv{(i_t^*-i_t)}{t}+\lambda_{5,t}\dv{(e_t^*-e_t)}{t}\\
    &=<\lambda_t,\dot{x}_t^*-\dot{x}_t>
\end{align*}
Also notice that,
\[\dv{\big(e^{-\theta t} \lambda_t\big)(x_t^*-x_t)}{t} = e^{-\theta t}<\lambda_t,\dot{x}_t^*-\dot{x}_t>+e^{-\theta t}<\lambda_t-\theta \lambda_t,x^*_t-x_t>\]
Therefore the second and third term of the Equation \ref{int1} can be simplify as,
\begin{align*}
&\int_{0}^{\infty} e^{-\theta t} [H(x_t^*,z_t,\lambda_t)-H(x_t,z_t,\lambda_t)]dt + \int_0^{\infty} e^{-\theta t}<\dot{\lambda}_t-\theta \lambda_t,x_t^*-x_t>dt \\ &=\int_0^\infty \bigg(\dv{\big(e^{-\theta t} \lambda_t\big)(x_t^*-x_t)}{t}\bigg) dt\\
&= lim_{t\to \infty} e^{-\theta t} \lambda_t(x_t^*-x_t) - \lambda_0 (x_0^*-x_0) \\
&=lim_{t\to \infty} e^{-\theta t} \lambda_t(x_t^*-x_t)\\
&=lim_{t\to \infty} e^{-\theta t} \lambda_{1,t}(k_t^*-k_t)+lim_{t\to \infty} e^{-\theta t} \lambda_{2,t}(h_t^*-h_t)+lim_{t\to \infty} e^{-\theta t} \lambda_{3,t}(s_t^*-s_t) \\ &+lim_{t\to \infty} e^{-\theta t} \lambda_{4,t}(i_t^*-i_t)+
lim_{t\to \infty} e^{-\theta t} \lambda_{5,t}(e_t^*-e_t)
\end{align*}
where $x_0^*=x_0$ from the initial conditions. 

\noindent Now we want to claim $lim_{t\to \infty} e^{-\theta t} \lambda_t(x_t^*-x_t)=0$. Using assumption (\ref{assboundk}) and similar argument given in proposition 4 in \cite{goenka2014infectious}, we conclude that \begin{equation}
   lim_{t\to \infty}e^{-\theta t} \lambda_{1,t}=lim_{t\to \infty}e^{-\theta t} \lambda_{2,t}=lim_{t\to \infty}e^{-\theta t} \lambda_{3,t}=lim_{t\to \infty}e^{-\theta t} \lambda_{4,t}=0 \label{trans1} 
\end{equation} Hence, 

\begin{align}
lim_{t\to \infty}e^{-\theta t} \lambda_{1,t}(k_t^*-k_t) &=lim_{t\to \infty}e^{-\theta t} \lambda_{2,t}(h_t^*-h_t) \nonumber\\
&=lim_{t\to \infty}e^{-\theta t} \lambda_{3,t}(s_t^*-s_t)\nonumber\\
&=lim_{t\to \infty}e^{-\theta t} \lambda_{4,t}(i_t^*-i_t)\nonumber\\
&=0 
\end{align}
From first order necessary conditions, 
\begin{equation}
    -\lambda_1-\lambda_3 \beta_1^* s^* i^* +\lambda_4 (\beta_1^* s^* i^*-\gamma_1^* i^*)+\lambda_5 E_1^*=0
\end{equation}
Since $E_1>0$ by assumptions and using Equation (\ref{trans1}), we conclude that \[lim_{t\to \infty}e^{-\theta t} \lambda_{5,t}=0\]
Thus,
\[lim_{t\to \infty}e^{-\theta t} \lambda_{5,t}(e_t^*-e_t)=0\]
Therefore,

\begin{equation}\label{int3}
    \int_{0}^{\infty} e^{-\theta t} [H(x_t^*,z_t,\lambda_t)-H(x_t,z_t,\lambda_t)]dt + \int_0^{\infty} e^{-\theta t}<\dot{\lambda}_t-\theta \lambda_t,x_t^*-x_t>dt=0
\end{equation}

Now using the Equations (\ref{int2}) and (\ref{int3}) we conclude that Equation (\ref{int1}) is non-negative. That is,
\begin{align*}
    &\int_{0}^{\infty} e^{-\theta t} [ H(x_t^*,z_t^*,\lambda_t)-H(x_t,z_t,\lambda_t)+<\dot{\lambda}_t-\theta
    \lambda_t,x_t^*-x_t>]dt \\
    &= \int_0^{\infty} e^{-\theta t} \big(u(c_t^*)-u(c_t)\big)dt\\
    &\ge 0
\end{align*}
Hence the result.
\end{proof}
Since endemic steady state solutions with positive health expenditure and positive direct investment to control the epidemic satisfy first order necessary conditions, we have following result.
\begin{corollary}
The endemic steady state with positive health expenditure and positive direct investment to control the epidemic is locally optimal.
\end{corollary}

\begin{corollary}
Given endemic state, let $x_t^*=(k_t^*,h_t^*,s_t^*,i_t^*,e_t^*), z_t^*=(c_t^*,m_t^*,A_t^*), $ and $\lambda_t=(\lambda_{1,t},\cdots,\lambda_{5,t})$ denote a path with positive health expenditure but zero direct investment to control the epidemic which satisfies necessary condition given in Equations (\ref{foncc}-\ref{focend}) along with the transversality conditions given in (\ref{foctrans1}-\ref{foctrans2}) and initial condition $x_0=(k_0,h_0,s_0,i_0,e_0)$. Then this path is a locally optimal for the social planner problem given in (\ref{social})
\end{corollary}
\begin{proof}
omitted
\end{proof}
\begin{corollary}
Given endemic state, let $x_t^*=(k_t^*,h_t^*,s_t^*,i_t^*,e_t^*), z_t^*=(c_t^*,m_t^*,A_t^*), $ and $\lambda_t=(\lambda_{1,t},\cdots,\lambda_{5,t})$ denote a path with zero health expenditure and positive direct investment to control the epidemic which satisfies necessary condition given in Equations (\ref{foncc}-\ref{focend}) along with the transversality conditions given in (\ref{foctrans1}-\ref{foctrans2}) and initial condition $x_0=(k_0,h_0,s_0,i_0,e_0)$. Then this path is a locally optimal for the social planner problem given in (\ref{social})
\end{corollary}
\begin{proof}
omitted
\end{proof}
\section{Numerical Results}
In this section, we present numerical simulations of the model in Python. To implement our model, we specify simple standard production function forms that meet the assumptions of the model, the transmission rate function form and the recovery rate function form.
We choose a Cobb-Douglas form for the production function of physical goods:
\begin{equation}\label{production}
f(k,l)=  k^{\psi} l^{1-\psi},
\end{equation}
where $0<\psi<1$. For healthy production function $g$, we choose 
\begin{equation}\label{health}
g(m)=\psi_3(m+\psi_1)^{\psi_2}-\psi_4\psi_1^{\psi_2}, 
\end{equation}
where $0<\psi_2<1,\psi_1\geq 0, \psi_3\geq 0$\cite{goenka2014infectious}. 

Our choice of transmission rate function form depends on  accumulative knowledge of the disease only:
\begin{equation}\label{trans11}
\beta(A,e,h)=\beta_0 e^{-\eta e}.
\end{equation}
Our second choice depends only health capital:
\begin{equation}\label{trans12}
\beta(A,e,h)=\beta_1+\beta_0 e^{-\eta h},
\end{equation}
Our third choice, we used transmission rate function given in \cite{boris}
\begin{equation}\label{trans13}
\beta(A,e,h)=\beta_0 e^{-\eta A e h},
\end{equation}
where $\eta\geq 0$. Here $\bar{\beta}=\beta_0$ or $\beta_1+\beta_0$. If there is no learning by controlling effect of the disease then we set $\eta=0$, in which case the transmission rate is a constant.

 We assume that the recovery rate depend on accumulative knowledge of the disease, health capital and capital to control the disease. we study the following examples:
 \begin{equation}
     \gamma(A,e,h)=\gamma_1-\gamma_0 e^{-\eta_2 A e h}
 \end{equation}
 where $\eta_2\geq 0$. Here $\underline{\gamma}=\gamma_1-\gamma_0>0$. If there is no learning by controlling effect of the disease then we set $\eta_2=0$, in which case the recovery rate is a constant.
 
 \begin{equation}
     \gamma(A,e,h)=\gamma_1-\gamma_0 e^{-\eta_2 h}.
 \end{equation}
 
We analyze the case when the parameter values produce the existence of an endemic steady state ($\mathcal{R}_0>1$) and strictly positive investments. Some values of parameters, we get them from the literature. We then choose the birth rate to match the approximate average fertility rate in following countries: France, Germany, Italy, United Kingdom, and United states  during the period of the  outbreak, 5 percent (World Bank). Unless otherwise noted, our numerical analyses use the following parameter values: 
$\psi=0.3,\phi_1 = 0.2,\, \phi_2= 0.5\, \phi_3=0.5,\phi_4 = 0.5,\phi_5, \varphi= 0.36, \delta_E= 0.05,\delta_K=0.05$, $\gamma_0=1,\gamma_1=1.01,\gamma_2=\gamma_3=1$,$a_1=a_2=0.023,a_3=1,p=0.5$,$\beta_0=0.023, \beta_1=0.023$  and $b= 0.0482$.
\subsection{Steady state without health expenditure}
We compute steady state for given $b\in [0.005, 0.13]$. For given $b$,
\[\underbar{l}=1-(0.5 \frac{b}{0.01}-\frac{b}{0.013})\in [1.13,1.35],\] we get $A^*=e^*=m^*=0,l^*=\underbar{l}$. Let $\psi=0.3$,then $k^*=\underbar{k}=\left(\frac{0.3}{\theta+0.045+b}\right)^{1/.7}$.
\subsection{Effects of Varying the Discount Rate}
We numerically  vary the social planner’s discount rate and true to see the  affects of the endemic steady state of our model with learning by controlling versus the model without learning by controlling. The discount rate, $\theta$ reflects the patience of the social planner, with a higher discount rate reflecting a less patient planner.

In the first case, we consider the model without disease control. In this case, the epidemiological variables are invariant to changes in the discount rate. The economic variables behave just as expected in a standard growth model, with a higher discount rate leading to a lower steady-state health capital stock and physical capital stock and therefore lower output and consumption.

The second case, we consider the model with disease. In this case, the epidemiological variables vary significantly with respect to the discount rate. Since changes in the discount rate affect learning by controlling and thus the disease transmission rate. 

The result of this is that a less patient planner chooses higher total and susceptible populations and lower infected and recovered populations. In terms of the economic variables, the stock of disease controlling experience varies dramatically with the discount rate. A less patient planner devotes far more resources to disease control. Because of the positive effect on the working population, this increased disease control does not come at the expense of output or consumption. Rather, output and consumption are both increasing in the discount rate. A final noteworthy aspect of the model with disease control is that, as the discount rate changes, general health investment and capital do not move monotonically. As the discount rate increases, they both initially fall and then rise.  We found that this may not be a good assumption once we allow for endogenous changes in health investment and disease control.
\begin{figure}[H]
\centering
\begin{minipage}{.5\textwidth}
  \centering
  \includegraphics[width=.6\linewidth]{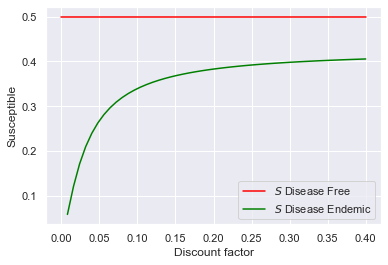}
  \caption{\footnotesize{Change in Susceptible Variable as Discount Date  Varies(Red Line - Disease-free Case; Green Line - Disease-endemic Case)}}
  \label{fig:S}
\end{minipage}%
\begin{minipage}{.5\textwidth}
  \centering
  \includegraphics[width=.6\linewidth]{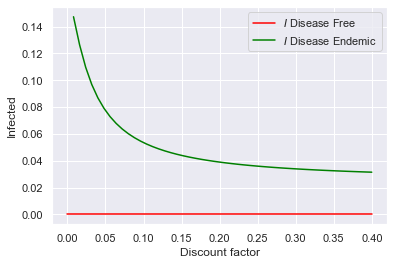}
  \caption{\footnotesize{Change in infected Variable as Discount Date  Varies(Red Line - Disease-free Case; Green Line - Disease-endemic Case)}}
  \label{fig:I}
\end{minipage}
\end{figure}

\begin{figure}[H]
\centering
\begin{minipage}{.5\textwidth}
  \centering
  \includegraphics[width=.6\linewidth]{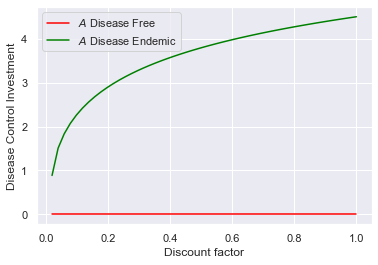}
  \caption{\footnotesize{Change in disease control investment Variable as Discount Date  Varies(Red Line - Disease-free Case; Green Line - Disease-endemic Case)}}
  \label{fig:A}
\end{minipage}%
\begin{minipage}{.5\textwidth}
  \centering
  \includegraphics[width=.5\linewidth]{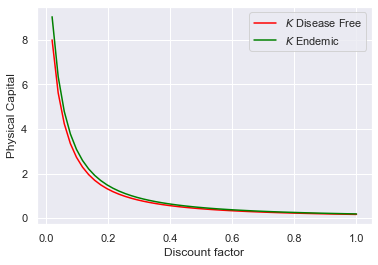}
  \caption{\footnotesize{Change in physical capital Variable as Discount Date  Varies(Red Line -Disease-free Case; Green Line - Disease-endemic Case)}}
  \label{fig:K}
\end{minipage}
\end{figure}

\begin{figure}[H]
\centering
\begin{minipage}{.5\textwidth}
  \centering
  \includegraphics[width=.5\linewidth]{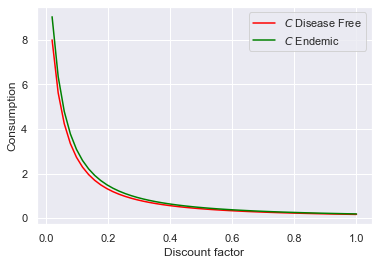}
  \caption{\footnotesize{Change in consumption Variable as Discount Date  Varies(Red Line - Disease-free Case; Green Line - Disease-endemic Case)}}
  \label{fig:C}
\end{minipage}%
\begin{minipage}{.5\textwidth}
  \centering
  \includegraphics[width=.5\linewidth]{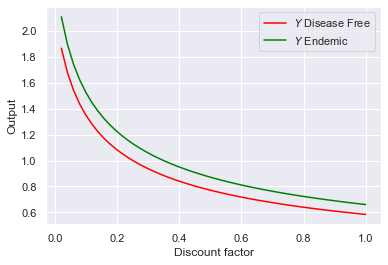}
  \caption{\footnotesize{Change in output Variable as Discount Date  Varies(Red Line - Disease-free Case; Green Line - Disease-endemic Case)}}
  \label{fig:O}
\end{minipage}
\end{figure}


\bibliographystyle{plain}
\bibliography{Ref}
\appendix
\end{document}